\documentclass[11pt]{amsproc}
\usepackage{amssymb,amsmath,amsthm,anysize}
\usepackage{longtable}
\usepackage{booktabs}
\usepackage{multirow,bigstrut,bigdelim,color}
\usepackage[shortlabels]{enumitem}
\usepackage{tabularx}
\usepackage[normalem]{ulem}
\usepackage{xcolor}

\marginsize{2cm}{2cm}{2cm}{2cm}

\newtheorem{Thm}{Theorem}[section]
\newtheorem{Cor}[Thm]{Corollary}

\newtheorem{Pro}[Thm]{Proposition}
\newtheorem{La}[Thm]{Lemma}
\newtheorem{Rem}[Thm]{Remark}
\newtheorem{Conj}[Thm]{Conjecture}
\newtheorem{hypothesis}[Thm]{Hypothesis}

\newtheorem{Prob}[Thm]{Problem}

\newenvironment{Prf}{\noindent\textbf{Proof.}}{\hfill $\Box$ \medskip}

\newcommand{\F}{\mathbb{F}}

\newcommand{\GL}{\mathrm{GL}}
\newcommand{\SL}{\mathrm{SL}}
\newcommand{\PGL}{\mathrm{PGL}}

\newcommand{\PSL}{\mathrm{PSL}}

\newcommand{\PSU}{\mathrm{PSU}}
\newcommand{\SO}{\mathrm{SO}}

\newcommand{\PGO}{\mathrm{PGO}}
\newcommand{\Sp}{\mathrm{Sp}}

\newcommand{\PSp}{\mathrm{PSp}}
\newcommand{\PGSp}{\mathrm{PGSp}}

\newcommand{\Aut}{\operatorname{Aut}}

\newcommand{\Alt}{\mathrm{Alt}}
\newcommand{\Sym}{\mathrm{Sym}}
\newcommand{\A}{\mathrm{A}}
\newcommand{\B}{\mathrm{B}}
\newcommand{\C}{\mathrm{C}}
\newcommand{\D}{\mathrm{D}}
\newcommand{\E}{\mathrm{E}}
\newcommand{\FF}{\mathrm{F}}
\newcommand{\G}{\mathrm{G}}
\newcommand{\inc}{\,\mathrm{I}\,}

\renewcommand\le{\leqslant}
\renewcommand\ge{\geqslant}

\setlist[description]{leftmargin=1ex,font=\normalfont\bfseries\space,style=nextline}

\begin{document}

\title[Simple groups, product actions, and generalised quadrangles]{Simple groups, product actions, and generalised quadrangles}
\author{John Bamberg, Tomasz Popiel, Cheryl E. Praeger}

\address{
Centre for the Mathematics of Symmetry and Computation\\
School of Mathematics and Statistics\\
The University of Western Australia\\
35 Stirling Highway, Crawley, W.A. 6009, Australia.\newline
Email: \texttt{\{john.bamberg, tomasz.popiel$^\dag$, cheryl.praeger$^\ddag$\}@uwa.edu.au} 
\newline $^\dag$ Current address: School of Mathematical Sciences, Queen Mary University of London, Mile End Road, London E1 4NS, United Kingdom. 
\newline $^\ddag$ Also affiliated with King Abdulaziz University, Jeddah, Saudi Arabia.
}

\thanks{The first author acknowledges the support of the Australian Research Council (ARC) Future Fellowship FT120100036. 
The second author acknowledges the support of the ARC Discovery Grant DP140100416. 
The research reported in the paper forms part of the ARC Discovery Grant DP140100416 of the third author. 
We thank Elisa Covato and Tim Burness for making available to us the results quoted in Section~\ref{ss:PA}, and Luke Morgan for helpful discussions.}



\maketitle

\begin{abstract}
The classification of flag-transitive generalised quadrangles is a long-standing open problem at the interface of finite geometry and permutation group theory. 
Given that all known flag-transitive generalised quadrangles are also point-primitive (up to point--line duality), it is likewise natural to seek a classification of the point-primitive examples. 
Working towards this aim, we are led to investigate generalised quadrangles that admit a collineation group $G$ preserving a Cartesian product decomposition of the set of points. 
It is shown that, under a generic assumption on $G$, the number of factors of such a Cartesian product can be at most four. 
This result is then used to treat various types of primitive and quasiprimitive point actions. 
In particular, it is shown that $G$ cannot have {\em holomorph compound} O'Nan--Scott type. 
Our arguments also pose purely group-theoretic questions about conjugacy classes in non-Abelian finite simple groups, and about fixities of primitive permutation groups. 
\end{abstract}

\section{Introduction} \label{ss:intro}

Generalised quadrangles are point--line incidence geometries introduced by Tits~\cite{Tits:1959cl} in an attempt to find geometric models for simple groups of Lie type. 
The {\em classical} generalised quadrangles arise in this way \cite[Section~3]{FGQ}. 
Each admits one of the simple classical groups $T = \text{PSp}(4,q) \cong \Omega_5(q)$, $\text{PSU}(4,q) \cong \text{P}\Omega_6^-(q)$ or $\text{PSU}(5,q)$ acting transitively on flags (incident point--line pairs). 
Moreover, the point and line stabilisers are certain maximal subgroups of $T$, so $T$ acts primitively on both points and lines. 
The classification of flag-transitive generalised quadrangles is a long-standing open problem. 
In addition to the classical families, only two other flag-transitive examples are known (up to point--line duality), each admitting an affine group acting point-primitively but line-imprimitively. 
Hence, all of the known flag-transitive generalised quadrangles are also point-primitive (up to duality), and so it is natural to seek a classification of the point-primitive examples. 
Indeed, this is arguably a more difficult problem, because one begins with essentially no information about the action of the collineation group on lines, nor any notion of what `incidence' means, whereas in a flag-transitive point--line geometry, points and lines correspond to cosets of certain subgroups of the collineation group, and incidence is determined by non-empty intersection of these cosets. 

Here we prove the following theorem. 
The abbreviations HS (holomorph simple), HC (holomorph compound), SD (simple diagonal), CD (compound diagonal), PA (product action), AS (almost simple) and TW (twisted wreath) refer to the possible types of non-affine primitive permutation group actions, in the sense of the O'Nan--Scott Theorem as stated in~\cite[Section~6]{CherylONS}. 
In the second column of Table~\ref{tab:Primitive2}, $\text{soc}(G)$ denotes the {\em socle} of the group $G$, namely the subgroup generated by its minimal normal subgroups. 
By $\text{fix}_\Omega(h)$ we mean the number of elements fixed by a permutation $h$ of the set $\Omega$, and $\mathsf{Q}^-(5,2)$ is the unique generalised quadrangle of order~$(2,4)$. 

\begin{Thm} \label{thm:QP-GQs}
If $\mathcal{Q}$ is a thick finite generalised quadrangle with a non-affine collineation group $G$ that acts primitively on the point set $\mathcal{P}$ of $\mathcal{Q}$, then the action of $G$ on $\mathcal{P}$ does not have O'Nan--Scott type HC, and the conditions in Table~\ref{tab:Primitive2} hold for the remaining O'Nan--Scott types.
\end{Thm}

\begin{table}[!ht]
\begin{small}
\begin{center}
\begin{tabular}{lll}
\toprule
Type & $\text{soc}(G)$ & Necessary conditions \\
\midrule
HS & $T \times T$ & $T$ is of Lie type $\A_5^\varepsilon$, $\A_6^\varepsilon$, $\B_3$, $\C_2$, $\C_3$, $\D_4^\varepsilon$, $\D_5^\varepsilon$, $\D_6^\varepsilon$, $\E_6^\varepsilon$, $\E_7$ or $\FF_4$ \\
\midrule
SD & $T^k$ & $T$ is a sporadic simple group or $T \cong \Alt_n$ with $n \le 18$; or \\
&& $T$ is an exceptional Lie type group; or \\
&& $T$ has Lie type $\A_1$, $\A_n^\varepsilon$ with $2 \le n \le 8$; $\B_n$ or $\C_n$ with $2 \le n \le 4$; or $\D_n^\varepsilon$ with $4 \le n \le 8$ \\
\cmidrule{2-3}
CD & $(T^k)^r$ & $r=2$ and $T \cong \Alt_n$ with $n \le 9$; or \\
&& $r=2$ and $T$ is a sporadic simple group with $T \not \cong \text{Suz}$, $\text{Co}_2$, $\text{Fi}_{22}$, $\text{Fi}_{23}$, $\text{B}$ or $\text{M}$; or \\
&& $r=2$ and $T$ has Lie type $\A_1$, $\A_2^\varepsilon$, $\A_3^\varepsilon$, $\B_2$, ${}^2\B_2$, ${}^2\FF_4$, $\G_2$ or ${}^2\G_2$; or \\
\cmidrule{3-3}
&& $r=3$ and $T \cong \text{J}_1$ or $T$ is of Lie type $\A_1$ or ${}^2\B_2$ \\
\midrule
PA & $T^r$ & $r=2$ and $\text{fix}_\Omega(h) < |\Omega|^{3/5}$ for all $h \in H \setminus \{1\}$; or \\
\cmidrule{3-3}
&& $3 \le r \le 4$, $T$ is a group of Lie type, and $\text{fix}_\Omega(h) < |\Omega|^{1-r/5}$ for all $h \in H \setminus \{1\}$; or \\
&& $3 \le r \le 4$, $H=T \cong \Alt_p$ with point stabiliser $p.\frac{p-1}{2}$ for a prime $p \equiv 3 \pmod 4$ \\
\midrule
AS & $T$ & $\text{fix}_\mathcal{P}(g) < |\mathcal{P}|^{4/5}$ for all $g \in G \setminus \{1\}$; or \\
&& $\mathcal{Q} = \mathsf{Q}^-(5,2)$ with $T \cong \PSU_4(2)$ \\ 
\midrule
TW & $T^r$ & $\text{fix}_\mathcal{P}(g) < |\mathcal{P}|^{4/5}$ for all $g \in G \setminus \{1\}$ \\
\bottomrule
\end{tabular}
\end{center}
\caption{\small 
Conditions for Theorem~\ref{thm:QP-GQs}. 
Here $T$ is a non-Abelian finite simple group, $k \ge 2$ and $r \ge 2$. 
If $G$ acts primitively of type CD (respectively PA) on $\mathcal{P}$, then $G \le H \wr \Sym_r$ for some primitive group $H \le \Sym(\Omega)$ of type SD (respectively AS) with socle $T^k$ (respectively~$T$). 
}\label{tab:Primitive2}
\end{small}
\end{table}

Note also that, in the notation for finite simple groups of Lie type used in Table~\ref{tab:Primitive2} (and throughout the paper), $\varepsilon = \pm$ and $\A_n^+=\A_n$, $\A_n^-={}^2\A_n$, $\D_n^+=\D_n$, $\D_n^-={}^2\D_n$, $\E_6^+=\E_6$, $\E_6^-={}^2\E_6$.

Before we proceed, a remark is in order about the assumption in Theorem~\ref{thm:QP-GQs} that $G$ not be an affine group. 
If $G$ is affine, then the generalised quadrangle $\mathcal{Q}$ necessarily arises from a so-called {\em pseudo-hyperoval} in a projective space $\text{PG}(3n-1,q)$ with $q$ even \cite{MR2201385}.
In joint work with Glasby \cite{OurAffine}, we were able to classify the generalised quadrangles admitting an affine group that acts primitively on points and transitively on lines: they are precisely the two flag-transitive, point-primitive, line-imprimitive generalised quadrangles mentioned above.  
However, without the extra assumption of transitivity on lines, the problem is equivalent to the classification of the pseudo-hyperovals that have an irreducible stabiliser. 
As explained in \cite[Remark~1.3]{OurAffine}, this latter problem would appear to be extremely difficult, and possibly intractable. 
It also has a rather different flavour to the cases treated in the present paper, and so we do not consider it further here.

Let us now establish some definitions and notation, before discussing further. 
By a {\em point--line incidence geometry} we mean a triple $\Gamma = (\mathcal{P},\mathcal{L},\inc)$, where $\mathcal{P}$ and $\mathcal{L}$ are sets whose elements are called {\em points} and {\em lines}, respectively, and $\inc \subseteq \mathcal{P} \times \mathcal{L}$ is a symmetric binary relation called {\em incidence}. 
We write $\Gamma = (\mathcal{P},\mathcal{L})$ instead of $(\mathcal{P},\mathcal{L},\inc)$ when we do not need to refer to the incidence relation explicitly.
Two points (respectively lines) of $\Gamma$ are said to be {\em collinear} (respectively {\em concurrent}) if they are incident with a common line (respectively point).
A {\em collineation} of $\Gamma$ is a permutation of $\mathcal{P} \cup \mathcal{L}$ that preserves $\mathcal{P}$ and $\mathcal{L}$ setwise and preserves the incidence relation.
By a {\em collineation group} of $\Gamma$ we mean a subgroup of the group of all collineations of $\Gamma$, which is called the {\em full collineation group}. 

A {\em generalised quadrangle} is a point--line incidence geometry $\mathcal{Q} = (\mathcal{P},\mathcal{L})$ that satisfies the following two axioms: (i) two distinct points are incident with at most one common line, and (ii) given a point $P$ and a line $\ell$ not incident with $P$, there exists a unique point incident with $\ell$ that is collinear with $P$. 
The second axiom implies that every pair in $\mathcal{P} \cup \mathcal{L}$ is contained in an ordinary quadrangle, and that $\mathcal{Q}$ contains no triangles. 
All generalised quadrangles considered in this paper are assumed to be finite, in the sense that $\mathcal{P}$ and $\mathcal{L}$ are finite sets.
If every point is incident with at least three lines, and every line is incident with at least three points, then $\mathcal{Q}$ is said to be {\em thick}. 
In this case, there exist constants $s \ge 2$ and $t \ge 2$ such that every point is incident with exactly $t+1$ lines and every line is incident with exactly $s+1$ points \cite[Corollary~1.5.3]{HvM}. 
The pair $(s,t)$ is called the {\em order} of $\mathcal{Q}$. 
Observe also that there is a natural concept of point--line duality for generalised quadrangles: if $(\mathcal{P},\mathcal{L})$ is a generalised quadrangle, then so is $(\mathcal{L},\mathcal{P})$; and if $(\mathcal{P},\mathcal{L})$ has order $(s,t)$, then $(\mathcal{L},\mathcal{P})$ has order $(t,s)$.

Let us now discuss Theorem~\ref{thm:QP-GQs} further. 
The primitive permutation groups on a finite set $\Delta$ are classified into eight types according to the O'Nan--Scott Theorem as presented in \cite[Section~6]{CherylONS}. 
In 2012, Bamberg et~al.~\cite{Bamberg:2012yf} showed that if a thick finite generalised quadrangle admits a collineation group $G$ that acts primitively on both points and lines, then $G$ must be an almost simple (AS type) group. 
That is, $T \le G \le \Aut(T)$ for some non-Abelian finite simple group $T$. 
Given that there exist point-primitive generalised quadrangles that are line-transitive but line-imprimitive, our initial aim was to extend the result of \cite{Bamberg:2012yf} by relaxing the line-primitivity assumption to line-transitivity. 
In addition to handling the affine (HA type) case with Glasby~\cite{OurAffine}, we were also able to show that no such examples arise if the point action has type HS or HC \cite{OurHSHC}.

Theorem~\ref{thm:QP-GQs} significantly strengthens and expands upon the results of \cite{Bamberg:2012yf,OurHSHC}. 
The idea behind its proof begins with the following observations. 
A primitive group $G \le \Sym(\Delta)$ of O'Nan--Scott type HC, CD, PA or TW preserves a Cartesian product decomposition $\Delta = \Omega^r$, for some set $\Omega$ and some $r \ge 2$. 
Therefore, in studying point-primitive generalised quadrangles, we are led in particular to consider generalised quadrangles with collineation groups that preserve a Cartesian product decomposition of the point set. 
The following theorem shows that the number of factors of such a decomposition becomes severely restricted under a fairly generic assumption on the group. 
Here a {\em semiregular} permutation group action is one in which only the identity element fixes a point, and if $H_1,\ldots,H_r$ are permutation groups on sets $\Omega_1,\ldots,\Omega_r$, respectively, then  the {\em product action} of the direct product $\prod_{i=1}^r H_i$ on the Cartesian product $\prod_{i=1}^r \Omega_i$ is the action $(\omega_1,\ldots,\omega_r)^{(h_1,\ldots,h_r)} = (\omega_1^{h_1},\ldots,\omega_r^{h_r})$. 
We also recall that a permutation group is said to act {\em regularly} if it acts transitively and semiregularly.

\begin{Thm} \label{thmPA-generic}
Let $\Omega_1,\ldots,\Omega_r$ be finite sets with $2 \le |\Omega_1| \le \cdots \le |\Omega_r|$, where $r \ge 1$, and let $H_i \le \Sym(\Omega_i)$ for each $i \in \{1,\ldots,r\}$. 
Assume further that $H_1$ is non-trivial and that its action on $\Omega_1$ is not semiregular. 
Suppose that $N = \prod_{i=1}^r H_i$ is a collineation group of a thick finite generalised quadrangle $\mathcal{Q} = (\mathcal{P},\mathcal{L})$ of order not equal to $(2,4)$, such that $\mathcal{P} = \prod_{i=1}^r \Omega_i$ and $N$ has the product action on $\mathcal{P}$. 
Then $r \le 4$, and every non-identity element of $H_1$ fixes less than $|\Omega_1|^{1-r/5}$ points of $\Omega_1$.
\end{Thm}

The proof of Theorem~\ref{thmPA-generic} relies on the existence of a non-identity element $h_1$ of $H_1$ that fixes at least one point of $\Omega_1$. 
If $r \ge 2$, one can then construct a collineation $(h_1,1,\ldots,1) \in N$ of $\mathcal{Q}$ that fixes at least $\prod_{i=2}^r |\Omega_i|$ points of the Cartesian product $\mathcal{P} = \prod_{i=1}^r \Omega_i$. 
Theorem~\ref{thmPA-generic} is then deduced from the following result, which bounds the number of points fixed by a non-identity collineation of an arbitrary thick finite generalised quadrangle. 
The proofs of both theorems are given in Section~\ref{sec:results}. 

\begin{Thm} \label{subGQLemma-general}
Let $\theta$ be a non-identity collineation of a thick finite generalised quadrangle $\mathcal{Q} = (\mathcal{P},\mathcal{L})$. 
Then either $\theta$ fixes less than $|\mathcal{P}|^{4/5}$ points of $\mathcal{Q}$, or $\mathcal{Q}$ is the unique generalised quadrangle $\mathsf{Q}^-(5,2)$ of order $(2,4)$ and $\theta$ fixes exactly $15$ of the $27$ points of $\mathcal{Q}$. 
\end{Thm}

\begin{Rem}
\textnormal{
Theorem~\ref{subGQLemma-general} improves a particular case of a recent result of Babai on automorphism groups of strongly regular graphs \cite[Theorem~1.7]{Babai}.
If $\mathcal{Q}$ has order $(s,t)$ then its collinearity graph, namely the graph with vertex set $\mathcal{P}$ and two vertices adjacent if and only if they are collinear in $\mathcal{Q}$, is a strongly regular graph with parameters $v = |\mathcal{P}| = (s+1)(st+1)$, $k = s(t+1)$, $\lambda = s-1$ and $\mu = t+1$. 
Roughly speaking, we have $v \approx s^2t$ and $k \approx st$, so the condition $k \le n^{3/4}$ in assertion~(b) of \cite[Theorem~1.7]{Babai} becomes $t \le s^2$, which is just Higman's inequality for generalised quadrangles (see Lemma~\ref{lemmaBasics}(ii)). 
Babai's result, which applies far more generally to strongly regular graphs that are {\em non-trivial}, {\em non-graphic} and {\em non-geometric}, therefore implies that a non-identity collineation $\theta$ of $\mathcal{Q}$ can fix at most $O(|\mathcal{P}|^{7/8})$ points. 
Theorem~\ref{subGQLemma-general} sharpens the $7/8$ exponent in this bound to $4/5$ in the case of collinearity graphs of generalised quadrangles. 
(Note also that assertion~(a) of \cite[Theorem~1.7]{Babai} sharpens the $7/8$ exponent to $5/6$ when, roughly, $t \ge s$: the condition $k \ge n^{2/3}$ roughly translates to $t \ge s$, and the corresponding bound is $O(\sqrt{kn})$, with $\sqrt{kn} \approx s^{3/2}t \ge (s^2t)^{5/6} \approx |\mathcal{P}|^{5/6}$ when $t \ge s$.)
}
\end{Rem}

To aid our discussion, let us now state the following immediate corollary of Theorem~\ref{thmPA-generic}. 

\begin{Cor} \label{corForPrim}
Let $\Omega$ be a finite set with $|\Omega| \ge 2$, and suppose that $H \le \Sym(\Omega)$ is non-trivial and not semiregular.  
Suppose that $N = H^r$, $r \ge 1$, is a collineation group of a thick finite generalised quadrangle $\mathcal{Q} = (\mathcal{P},\mathcal{L})$ of order not equal to $(2,4)$, such that $\mathcal{P} = \Omega^r$ and $N$ has the product action on $\mathcal{P}$. 
Then $r \le 4$, and every non-identity element of $H$ fixes less than $|\Omega|^{1-r/5}$ points of $\Omega$.
\end{Cor}

We apply Corollary~\ref{corForPrim} to groups $N$ that arise as subgroups of certain types of primitive groups. 
This in turn motivates certain questions about non-Abelian finite simple groups. 
As illustration, consider the case where $\Omega = T$ for some non-Abelian finite simple group $T$, with $H=T \times T$ acting on $\Omega$ via $\omega^{(x,x')} = x^{-1}\omega x'$. 
This situation arises when $N$ is the socle (the subgroup generated by the minimal normal subgroups) of a primitive group of type HS ($r=1$) or HC ($r \ge 2$). 
If $x'=x$ then the element $(x,x') = (x,x) \in H$ fixes precisely $|C_T(x)|$ points of $\Omega$, where $C_T(x)$ is the centraliser of $x$ in $T$. 
Corollary~\ref{corForPrim} therefore implies that $r \le 4$, and that $|C_T(x)| < |T|^{1-r/5}$ for all $x \in T \setminus \{1\}$. 
We therefore ask which non-Abelian finite simple groups $T$ satisfy this condition. 
If $r=4$ then we require that $|C_T(x)| < |T|^{1/5}$ for all $x \in T \setminus \{1\}$, which is false for every non-Abelian finite simple group $T$. 
Indeed, it is well known that every non-Abelian finite simple group $T$ contains an involution $x$ with $|C_T(x)| > |T|^{1/3}$ (in fact, {\em every} involution in $T$ has this property \cite[Proposition~2.4]{LiebeckShalevFixity}). 
For $r \in \{1,2,3\}$, we verify the following result in Section~\ref{sec:groups}. 
Although this result follows from routine calculations, we include it here in case it proves to be a convenient reference.

\begin{table}[!t]
\begin{small}
\begin{center}
\begin{tabular}{llll}
\toprule
$T$ & $r=1$ & $r=2$ & $r=3$ \\
\midrule
$\Alt_n$ & $5 \le n \le 18$ & $5 \le n \le 9$ & $5 \le n \le 6$ \\
sporadic & any & any except $\text{Suz}$, $\text{Co}_2$, $\text{Fi}_{22}$, $\text{Fi}_{23}$, $\text{B}$, $\text{M}$ & $\text{J}_1$ \\
exceptional Lie type & any & ${}^2\FF_4(q)$, $\G_2(q)$, ${}^2\G_2(q)$, ${}^2\B_2(q)$ & ${}^2\B_2(q)$ \\
$\PSL_{n+1}(q)$ & $1 \le n \le 8$ & $1 \le n \le 3$ & $n=1$, $q \neq 7$ \\
$\PSU_{n+1}(q)$ & $2 \le n \le 8$ & $2 \le n \le 3$ & --- \\
$\PSp_{2n}(q)$ or $\Omega_{2n+1}(q)$ & $2 \le n \le 4$ & $n=2$ & --- \\
$\text{P}\Omega_{2n}^\pm(q)$ & $4 \le n \le 8$ & --- & --- \\
\bottomrule
\end{tabular}
\end{center}
\caption{\small Possibilities for a non-Abelian finite simple group $T$ with the property that $|C_T(x)| < |T|^{1-r/5}$ for all $x \in T \setminus \{1\}$, for $r$ equal to one of $1$, $2$ or $3$.}\label{tab:summary4}
\end{small}
\end{table}

\begin{Pro} \label{classes}
Let $r \in \{1,2,3\}$ and let $T$ be a non-Abelian finite simple group. 
Then either $|C_T(x)| \ge |T|^{1-r/5}$ for some $x \in T \setminus \{1\}$, or $T$ is one of the groups listed in Table~\ref{tab:summary4}.
\end{Pro}

Our new results about generalised quadrangles with point-primitive collineation groups are proved in Sections~\ref{sec:SDCD}--\ref{ss:PA}. 
Corollary~\ref{corForPrim} is applied not only to actions of type HS or HC as illustrated above, but also to types SD, CD and PA. 
In particular, the proof of Theorem~\ref{thm:QP-GQs} is free from the Classification of Finite Simple Groups (CFSG) to the extent that, for $G$ of type HC, CD or PA with socle $T^r \times T^r$, $T^{(k-1)r}$ or $T^r$ respectively, the proof that $r \le 4$ depends only on Corollary~\ref{corForPrim}. 
The CFSG is, however, needed to prove Proposition~\ref{classes} and some of the results in Section~\ref{sec:PA}. 
For type PA, the group $H$ in Corollary~\ref{corForPrim} is an almost simple primitive group, so we are led to consider lower bounds on the {\em fixity} of such a group, namely the maximum number of fixed points of a non-identity element. 
In Section~\ref{ss:PA}, we discuss how refinements of a recent result of Liebeck and Shalev~\cite[Theorem~4]{LiebeckShalevFixity} on this problem, currently being carried out by Elisa Covato at the University of Bristol as part of her PhD research \cite{ElisaPhD}, can be adapted to further improve the bound $r \le 4$ in this case. 
In particular, for $r \in \{3,4\}$ we are able to show that $T$ cannot be a sporadic simple group, and to rule out the case $T \cong \Alt_n$ except in one specific action when $n$ is a prime congruent to $3$~modulo~$4$ (see Table~\ref{tab:Primitive2}). 
The proof of Theorem~\ref{thm:QP-GQs} is presented in Section~\ref{ss:mainProof}. 

Section~\ref{ss:discussion} concludes the paper with a discussion and some open problems. 
In light of the growing body of work towards a classification of point-primitive generalised quadrangles, and the possible avenues outlined in Remark~\ref{HSremark}, Remark~\ref{RemPA} and Section~\ref{ss:discussion} for attacking the cases left open by Theorem~\ref{thm:QP-GQs}, we feel that the following conjecture can be made with a reasonable amount of confidence.

\begin{Conj}
If a thick finite generalised quadrangle $\mathcal{Q}$ admits a collineation group $G$ that acts primitively on the point set of $\mathcal{Q}$, then $G$ is either affine or almost simple.
\end{Conj}

\section{Bounding the number of points fixed by a collineation} \label{sec:results}

The facts summarised in the following lemma are well known. 
(The existence of an order is proved in \cite[Corollary~1.5.3]{HvM}, and proofs of assertions (i)--(iii) may be found in \cite[Section~1.2]{FGQ}.)

\begin{La} \label{lemmaBasics}
Let $\mathcal{Q}$ be a thick finite generalised quadrangle. 
Then $\mathcal{Q}$ has an order $(s,t)$, and the following properties hold:
\begin{itemize}
\item[\textnormal{(i)}] $\mathcal{Q}$ has $(s+1)(st+1)$ points and $(t+1)(st+1)$ lines,
\item[\textnormal{(ii)}] $s^{1/2} \le t \le s^2 \le t^4$ (Higman's inequality),
\item[\textnormal{(iii)}] $s+t$ divides $st(st+1)$.
\end{itemize}
\end{La}

A point--line incidence geometry $\mathcal{S} = (\mathcal{P},\mathcal{L},\inc)$ is called a {\em grid} if there exist positive integers $s_1$ and $s_2$ such that $\mathcal{P}$ can be written as $\{ P_{ij} \mid 0 \le i \le s_1, 0 \le j \le s_2 \}$, $\mathcal{L}$ can be written as $\{ \ell_k \mid 0 \le k \le s_1 \} \cup \{ \ell_k' \mid 0 \le k \le s_2 \}$, and we have $P_{ij} \inc \ell_k$ if and only if $i=k$, and $P_{ij} \inc \ell_k'$ if and only if $j=k$.  
Each point of $\mathcal{S}$ is then incident with exactly two lines, and $|\mathcal{P}| = (s_1+1)(s_2+1)$. 
Let us say that such a grid has {\em parameters} $s_1$ and $s_2$. 
Note that a grid with parameters $s_1=s_2$ is a generalised quadrangle of order $(s_1,1)$. 
A {\em dual grid} is defined analogously, by swapping the roles of points and lines. 
That is, there exist positive integers $t_1$ and $t_2$ such that $\mathcal{L}$ can be written as $\{ \ell_{ij} \mid 0 \le i \le t_1, 0 \le j \le t_2 \}$, $\mathcal{P}$ can be written as $\{ P_i \mid 0 \le i \le t_1 \} \cup \{ P_j' \mid 0 \le j \le t_2 \}$, $P_k \inc \ell_{ij}$ if and only if $i=k$, and $P_k' \inc \ell_{ij}$ if and only if $j=k$.  
In this case, each line is incident with exactly two points, and $|\mathcal{P}| = (t_1+1)+(t_2+1)$. 
Let us say that such a dual grid has {\em parameters} $t_1$ and $t_2$.

If $\theta$ is a collineation of a generalised quadrangle $\mathcal{Q} = (\mathcal{P},\mathcal{L})$, then it makes sense to consider the point--line incidence geometry $\mathcal{Q}_\theta = (\mathcal{P}_\theta,\mathcal{L}_\theta)$ with $\mathcal{P}_\theta = \{ P \in \mathcal{P} \mid P^\theta = P \}$, $\mathcal{L}_\theta = \{ \ell \in \mathcal{L} \mid \ell^\theta = \ell \}$, and incidence inherited from $\mathcal{Q}$. 
Here we call $\mathcal{Q}_\theta$ the {\em substructure of $\mathcal{Q}$ fixed by $\theta$}. 
It may happen that $\mathcal{Q}_\theta$ is a grid or a dual grid, or a generalised quadrangle. 
More specifically, we have the following result, based on the description of the possible structures of $\mathcal{Q}_\theta$ given by Payne and Thas \cite[2.4.1]{FGQ}.

\begin{La} \label{substructureLemma}
Let $\mathcal{Q} = (\mathcal{P},\mathcal{L})$ be a thick finite generalised quadrangle of order $(s,t)$. 
Let $\theta$ be a non-identity collineation of $\mathcal{Q}$, and let $\mathcal{Q}_\theta = (\mathcal{P}_\theta,\mathcal{L}_\theta)$ be the substructure of $\mathcal{Q}$ fixed by $\theta$. 
Then at least one of the following conditions holds. 
\begin{itemize}
\item[\textnormal{(i)}] $\mathcal{P}_\theta$ is empty. 
\item[\textnormal{(ii)}] $\mathcal{L}_\theta$ is empty and $\mathcal{P}_\theta$ is a set of pairwise non-collinear points. 
In particular, $|\mathcal{P}_\theta| \le st+1$.
\item[\textnormal{(iii)}] All points of $\mathcal{Q}_\theta$ are incident with a common line, and $|\mathcal{P}_\theta| \le s+1$.
\item[\textnormal{(iv)}] All points of $\mathcal{Q}_\theta$ are collinear with a common point, and $|\mathcal{P}_\theta| \le s(t+1)+1$. 
\item[\textnormal{(v)}] $\mathcal{Q}_\theta$ is a grid. In this case, either $|\mathcal{P}_\theta| = (s+1)^2$ and $s \le t$, or $|\mathcal{P}_\theta| < s^2$.
\item[\textnormal{(vi)}] $\mathcal{Q}_\theta$ is a dual grid, and $|\mathcal{P}_\theta| \le 2(t+1)$.
\item[\textnormal{(vii)}] $\mathcal{Q}_\theta$ is a thick generalised quadrangle, and $|\mathcal{P}_\theta| \le (s+1)(t+1)$. 
\end{itemize}
In particular, either $|\mathcal{P}_\theta| \le (s+1)(t+1)$; or $s \ge t+3$, $\mathcal{Q}_\theta$ is a grid and $|\mathcal{P}_\theta| < s^2$.
\end{La}

\begin{Prf}
The possible structures (i)--(vii) of $\mathcal{Q}_\theta$ are given by \cite[2.4.1]{FGQ}. 
We verify the claimed upper bounds for $|\mathcal{P}_\theta|$. 
The bounds in cases~(iii) and~(iv) are immediate, because every line of $\mathcal{Q}$ is incident with exactly $s+1$ points, and every point of $\mathcal{Q}$ is incident with exactly $t+1$ lines. 
For case~(ii), note \cite[Section~2.7]{FGQ} that the maximum size of a set of pairwise non-collinear points in $\mathcal{Q}$ is $st+1$. 
For (v), if $\mathcal{Q}_\theta$ is a dual grid with parameters $t_1$ and $t_2$, then $t_1 \le t$ and $t_2 \le t$, and hence $|\mathcal{P}_\theta| \le 2(t+1)$. 

Now suppose that $\mathcal{Q}_\theta$ is a grid with parameters $s_1$ and $s_2$, noting that $s_1 \le s$ and $s_2 \le s$, and assuming (without loss of generality) that $s_1 \ge s_2$. 
If $s_1=s_2=s$ then $|\mathcal{P}_\theta| = (s+1)^2$, and $\mathcal{Q}_\theta$ is a generalised quadrangle of order $(s,1)$, so \cite[2.2.2(i)]{FGQ} implies that $s \le t$. 
The case $s_2 = s-1$ cannot occur, because if $\theta$ fixes $s$ points incident with a line then it must also fix the final point; and if $s_2 \le s-2$ then $|\mathcal{P}_\theta| \le (s+1)(s-1) < s^2$.  
Finally, suppose that $\mathcal{Q}$ is a thick finite generalised quadrangle, and let $(s',t')$ denote its order. 
Then $|\mathcal{P}_\theta| = (s'+1)(s't'+1)$ by Lemma~\ref{lemmaBasics}(i).
If $t'=t$ then $s'<s$ because $\theta \neq 1$, so \cite[2.2.1]{FGQ} implies that $s't = s't' \le s$, and hence $|\mathcal{P}_\theta| \le (s/t+1)(s+1) \le (t^2/t+1)(s+1) = (s+1)(t+1)$, where for the second inequality we use Lemma~\ref{lemmaBasics}(ii). 
If $t' < t$ then the dual statement of \cite[2.2.1]{FGQ} yields $s't' \le t$, so $|\mathcal{P}_\theta| = (s'+1)(s't'+1) \le (s+1)(t+1)$. 

The final assertion is deduced by comparing the upper bounds on $|\mathcal{P}_\theta|$ established in each case. 
We observe that $|\mathcal{P}_\theta| \le (s+1)(t+1)$ except possibly in the second case of~(v), where our bound is $|\mathcal{P}_\theta| < s^2$. 
However, if $s \le t+2$ then in this case we have $|\mathcal{P}_\theta| < s^2 < (s+1)(t+1)$. 
\end{Prf}

\begin{Rem}
\textnormal{
We mention a paper of Frohardt and Magaard \cite[Section~1.3]{FrohardtMagaard}, in which results analogous to Lemma~\ref{substructureLemma} are obtained for generalised $d$-gons with $d \in \{6,8\}$ (that is, generalised hexagons and octagons). 
The known examples of such geometries admit point- and line-primitive actions of almost simple groups with socle ${}^3\D_4(q)$ or $\G_2(q)$ (for $d=6$) and ${}^2\FF_4(q)$ (for $d=8$). 
Frohardt and Magaard use the aforementioned results to determine upper bounds for fixities of primitive actions of groups $G$ with generalised Fitting subgroup ${}^3\D_4(q)$, $\G_2(q)$ or ${}^2\FF_4(q)$ (and they also treat the other exceptional Lie type groups of Lie rank $1$ or $2$). 
By comparison, we instead apply Lemma~\ref{substructureLemma} to determine which groups might act primitively on the points of a generalised $d$-gon (with $d=4$ in our case). 
(We remark that we have also investigated point-primitive generalised hexagons and octagons, although via different methods than in the present paper \cite{HexOct,HexOctLuke}.)
}
\end{Rem}

We now use Lemma~\ref{substructureLemma} to prove Theorem~\ref{subGQLemma-general}, from which we deduce Theorem~\ref{thmPA-generic}. 
\medskip

\noindent {\bf Proof of Theorem~\ref{subGQLemma-general}}
Let $(s,t)$ be the order of $\mathcal{Q}$, and let $\mathcal{Q}_\theta = (\mathcal{P}_\theta,\mathcal{L}_\theta)$ be the substructure of $\mathcal{Q}$ fixed by $\theta$. 
We must show that either $|\mathcal{P}_\theta| < |\mathcal{P}|^{4/5}$, or $(s,t)=(2,4)$ and $|\mathcal{P}_\theta| = 15$ for $\mathcal{Q} = \mathsf{Q}^-(5,2)$.

First suppose that $s \neq 2$. 
By Lemma~\ref{substructureLemma}, we have either $|\mathcal{P}_\theta| < s^2$ or $|\mathcal{P}_\theta| \le (s+1)(t+1)$. 
If $|\mathcal{P}_\theta| < s^2$ then $|\mathcal{P}_\theta| < |\mathcal{P}|^{4/5}$ since $|\mathcal{P}| = (s+1)(st+1) > s^2t \ge s^{5/2}$ by Lemma~\ref{lemmaBasics}. 
If $|\mathcal{P}_\theta| \le (s+1)(t+1)$ then it suffices to show that the function $f(s,t) = ((s+1)(st+1))^{4/5} - (s+1)(t+1)$ is positive for all $s \ge 3$, for all $s^{1/2} \le t \le s^2$. 
This is readily checked when $s\in\{3,4\}$, so assume that $s \ge 5$. 
Thinking of $s$ and $t$ as real variables, we have 
\[
\frac{\partial f}{\partial t}(s,t) = \frac{(s+1)(4s-h(s,t))}{h(s,t)}, \quad \text{where} \quad h(s,t) = 5((s+1)(st+1))^{1/5}.
\]
Since $s$ and $t$ are positive, this derivative is positive if and only if $4s - h(s,t) > 0$.
Since $s \ge 5$ and $2 \le t \le s^2$, we have $h(s,t) \le 5 (\tfrac{6}{5} s)^{1/5} (\tfrac{11}{10} st)^{1/5} \le 5 (\tfrac{33}{25})^{1/5} s^{4/5}$. 
Hence, $4s - h(s,t) \ge s^{4/5} (4s^{1/5} - 5(\tfrac{33}{25})^{1/5})$. 
The right-hand side of this inequality is positive if $s > (\tfrac{5}{4})^5(\tfrac{33}{25}) = \tfrac{4125}{1024} \approx 4.028$, and so certainly $\tfrac{\partial f}{\partial t}(s,t) > 0$ when $s \ge 5$ and $s^{1/2} \le t \le s^2$. 
Since $f(s,t) \ge f(s,s^{1/2})$ and $f(s,s^{1/2}) > 0$ for $s \ge 5$, it follows that $f(s,t) > 0$ for all $s \ge 5$, for all $s^{1/2} \le t \le s^2$.

Now suppose that $s=2$. 
Then $t \in \{2,4\}$ by Lemma~\ref{lemmaBasics}.
There exist unique generalised quadrangles of orders $(2,2)$ and $(2,4)$, namely the symplectic space $\mathsf{W}(3,2)$ and the elliptic quadric $\mathsf{Q}^-(5,2)$, respectively \cite[5.2.3 and 5.3.2]{FGQ}. 
The full collineation groups of these generalised quadrangles are $\text{P}\Gamma\text{Sp}_4(2)$ and $\text{P}\Gamma\text{U}_4(2)$, respectively. 
One may use the package {\sf FinInG}~\cite{FinInG} in the computer algebra system {\sf GAP}~\cite{GAP} to check that every non-identity collineation of $\mathsf{W}(3,2)$ fixes at most $7$ points. 
Since $\mathsf{W}(3,2)$ has a total of $15$ points and $15^{4/5} \approx 8.73 > 7$, the claimed inequality $|\mathcal{P}_\theta| < |\mathcal{P}|^{4/5}$ holds for every non-identity collineation $\theta$ in this case.
On the other hand, there exist $36$ non-identity collineations of $\mathsf{Q}^-(5,2)$ that fix $15$ points, but the total number of points of $\mathsf{Q}^-(5,2)$ is $27$ and $27^{4/5} \approx 13.97 < 15$. 
We also remark that the substructure fixed by such a collineation is, in fact, a generalised quadrangle of order $(2,2)$.
Every other non-identity collineation of $\mathsf{Q}^-(5,2)$ fixes at most $9$ points.
\hfill $\Box$ \medskip

\noindent {\bf Proof of Theorem~\ref{thmPA-generic}.}
Since the action of $H_1$ on $\Omega_1$ is not semiregular, there exists $h_1 \in H_1 \setminus \{1\}$ fixing at least one point of $\Omega_1$. 
Let $f_1$ be the number of points of $\Omega_1$ fixed by $h_1$.
Let $\theta = (h_1,1,\ldots,1) \in N$, and let $f$ be the number of points of $\mathcal{Q}$ fixed by $\theta$.
If $r=1$ then Theorem~\ref{subGQLemma-general} implies that $f_1 = f < |\mathcal{P}|^{4/5} = |\Omega_1|^{4/5}$. 
If $r \ge 2$ then Theorem~\ref{subGQLemma-general} gives $f_1 \prod_{i=2}^r |\Omega_i| = f < |\mathcal{P}|^{4/5} = (\prod_{i=1}^r |\Omega_i|)^{4/5}$, so $f_1 (\prod_{i=2}^r |\Omega_i|)^{1/5} < |\Omega_1|^{4/5}$. 
Since $(\prod_{i=2}^r |\Omega_i|)^{1/5} \ge |\Omega_1|^{(r-1)/5}$, it follows that $f_1 < |\Omega_1|^{1-r/5}$. 
In particular, $1-r/5 > 0$ because $f_1 \ge 1$, and so $r \le 4$.
\hfill $\Box$ \medskip

We also use Lemma~\ref{substructureLemma} to sharpen the $4/5$ exponent bound in Theorem~\ref{subGQLemma-general} in some special cases. 
The proofs are just modifications of the proof of Theorem~\ref{subGQLemma-general}, but since the details are somewhat tedious to check, we include them in Appendix~\ref{ss:app} to save the reader having to reproduce them. 
We also remark that the exponent $94/125 = 0.752$ in case (i) of Proposition~\ref{Prop.752} could be changed to $3/4+\epsilon$ for any $\epsilon > 0$ at the expense of increasing the upper bound on $s$ in case (ii), but that this would not have been useful for our arguments in Section~\ref{sec:PA}.

\begin{Pro} \label{Prop7/9}
Let $\mathcal{Q} = (\mathcal{P},\mathcal{L})$ be a finite generalised quadrangle of order $(s,t)$, let $\theta$ be any non-identity collineation of $\mathcal{Q}$, and let $\mathcal{Q}_\theta = (\mathcal{P}_\theta,\mathcal{L}_\theta)$ be the substructure of $\mathcal{Q}$ fixed by $\theta$. 
Then either
\begin{itemize}
\item[\textnormal{(i)}] $|\mathcal{P}_\theta| < |\mathcal{P}|^{7/9}$, 
\item[\textnormal{(ii)}] $s \in \{2,3\}$, $t=s^2$ and $\mathcal{Q}_\theta$ is a generalised quadrangle of order $(s,s)$, or 
\item[\textnormal{(iii)}] $s \ge t+3$, $\mathcal{Q}_\theta$ is a grid and $|\mathcal{P}_\theta| < s^2$.
\end{itemize}
\end{Pro}

\begin{Pro} \label{Prop.752}
Let $\mathcal{Q} = (\mathcal{P},\mathcal{L})$ be a finite generalised quadrangle of order $(s,t)$, let $\theta$ be any non-identity collineation of $\mathcal{Q}$, and let $\mathcal{Q}_\theta = (\mathcal{P}_\theta,\mathcal{L}_\theta)$ be the substructure of $\mathcal{Q}$ fixed by $\theta$. 
Then either
\begin{itemize}
\item[\textnormal{(i)}] $|\mathcal{P}_\theta| < |\mathcal{P}|^{94/125}$, 
\item[\textnormal{(ii)}] $s < 2.9701 \times 10^{15}$, or 
\item[\textnormal{(iii)}] $s \ge t+3$, $\mathcal{Q}_\theta$ is a grid and $|\mathcal{P}_\theta| < s^2$.
\end{itemize}
\end{Pro}

\begin{Pro} \label{t=s+2}
Let $\mathcal{Q} = (\mathcal{P},\mathcal{L})$ be a finite generalised quadrangle of order $(s,t)$, let $\theta$ be any non-identity collineation of $\mathcal{Q}$, and let $\mathcal{P}_\theta$ denote the set of points fixed by $\theta$.
Suppose that $t=s+2$. 
Then $|\mathcal{P}_\theta| < |\mathcal{P}|^{7/9}$ if $s \ge 3$, and $|\mathcal{P}_\theta| < |\mathcal{P}|^{13/18}$ if $s \ge 5$.
\end{Pro}

\section{Centraliser orders in non-Abelian finite simple groups} \label{sec:groups}

Here we verify a series of lemmas about centraliser orders in non-Abelian finite simple groups, from which Proposition~\ref{classes} is deduced. 
Specifically, we need to know which non-Abelian finite simple groups $T$ contain non-identity elements $x$ with `large' centralisers, in the sense that $|C_T(x)| > |T|^{1-r/5}$ for $r$ equal to one of $1$, $2$ or $3$. 
This question is readily and exactly answered for alternating groups and sporadic simple groups in the following two lemmas.
Note that we treat the Tits group ${}^2\FF_4(2)'$ in Lemma~\ref{CTgSpor} along with the sporadic groups.

\begin{La} \label{CTgAltn}
Let $T \cong \Alt_n$ with $n \ge 5$. 
Then
\begin{itemize}
\item[\textnormal{(i)}] $|C_T(x)| < |T|^{4/5}$ for all $x \in T \setminus \{1\}$ if and only if $n \le 18$,
\item[\textnormal{(ii)}] $|C_T(x)| < |T|^{3/5}$ for all $x \in T \setminus \{1\}$ if and only if $n \le 9$,
\item[\textnormal{(iii)}] $|C_T(x)| < |T|^{2/5}$ for all $x \in T \setminus \{1\}$ if and only if $n \le 6$.
\end{itemize}
\end{La}

\begin{Prf}
If $n \ge 19$ and $x \in T$ is a $3$-cycle, then $|C_T(x)| = \tfrac{3}{2}(n-3)! > (\tfrac{1}{2}n!)^{4/5} = |T|^{4/5}$. 
The remaining assertions are readily verified using \sf{GAP}~\cite{GAP}.
\end{Prf}

\begin{La} \label{CTgSpor}
Let $T$ be either a sporadic finite simple group or the Tits group ${}^2\FF_4(2)'$. 
Then
\begin{itemize}
\item[\textnormal{(i)}] $|C_T(x)| < |T|^{4/5}$ for all $x \in T \setminus \{1\}$, 
\item[\textnormal{(ii)}] $|C_T(x)| < |T|^{3/5}$ for all $x \in T \setminus \{1\}$ if and only if $T \not\cong \textnormal{Suz}$, $\textnormal{Co}_2$, $\textnormal{Fi}_{22}$, $\textnormal{Fi}_{23}$, $\textnormal{B}$ or $\textnormal{M}$,
\item[\textnormal{(iii)}] $|C_T(x)| < |T|^{2/5}$ for all $x \in T \setminus \{1\}$ if and only if $T \cong \textnormal{J}_1$.
\end{itemize}
\end{La}

\begin{Prf}
This is readily verified upon checking maximal centraliser orders in the ATLAS~\cite{ATLAS}. 
\end{Prf}

Next we consider the exceptional Lie type groups, namely those of type $\E_8$, $\E_7$, $\E_6^\varepsilon$ (where $\varepsilon = \pm$), $\FF_4$, ${}^2\FF_4$, $\G_2$, ${}^2\G_2$, ${}^3\D_4$ or ${}^2\B_2$. 
Note that we make no attempt to check the converse of assertion~(i) (although this could be done using standard references including those cited here).

\begin{La} \label{centralisersExceptional}
Let $T$ be a finite simple group of exceptional Lie type. 
\begin{itemize}
\item[\textnormal{(i)}] If $T$ has type $\E_8$, $\E_7$, $\E_6^\varepsilon$, $\FF_4$ or $^3\D_4$, then there exists $x \in T \setminus \{1\}$ with $|C_T(x)| > |T|^{3/5}$.
\item[\textnormal{(ii)}] $|C_T(x)| < |T|^{2/5}$ for all $x \in T \setminus \{1\}$ if and only if $T$ has type ${}^2\B_2$.
\end{itemize}
\end{La}

\begin{Prf}
(i) For $T \cong \E_8(q)$, $\FF_4(q)$ or ${}^3\D_4(q)$, take $x\in T$ to be a unipotent element of type $\A_1$ in the sense of \cite[Tables~22.2.1 and~22.2.4]{LiebeckSeitz} and \cite{Spalt3D4}, respectively. 
Then $|C_T(x)| = q^{57}|\E_7(q)|$, $q^{15}|\C_3(q)|$ or $q^{12}(q^6-1)$, respectively, and it is readily checked that $|C_T(x)| > |T|^{3/5}$ in each case. 
Now suppose that $T \cong \E_7(q)$ or $\E_6^\varepsilon(q)$, and write $G := \text{Inndiag}(T)$. 
Take $x\in T$ to be a unipotent element of type $\A_1$ in the sense of \cite[Tables~22.2.2 and~22.2.3]{LiebeckSeitz}, respectively. 
Then $x^T = x^G$ by \cite[Corollary~17.10]{BG}, so $|C_T(x)| = |C_G(x)|/|G:T| = q^{33}|\D_6(q)|/\gcd(2,q-1)$ or $q^{21}|\A_5^\varepsilon(q)|/\gcd(3,q-\varepsilon)$, respectively, and again one can check that $|C_T(x)| > |T|^{3/5}$ in each case. 

(ii) If $T \cong {}^2\B_2(q)$ then $|C_T(x)| \le q^2 < (q^2(q^2+1)(q-1))^{2/5} = |T|^{2/5}$ for all $x \in T \setminus \{1\}$ \cite{Suzuki2}. 
It remains to check that $|C_T(x)| > |T|^{2/5}$ for some $x \in T \setminus \{1\}$ when $T$ has type ${}^2\FF_4$, $\G_2$ or ${}^2\G_2$. 
In these respective cases, take $x$ to be a unipotent element of type $(\tilde\A_1)_2$, $\A_1$ or $(\tilde\A_1)_3$ in the sense of \cite[Tables~22.2.5--22.2.7]{LiebeckSeitz}, so that $|C_T(x)| = q^{10}|{}^2\B_2(q)|$, $q^5|\A_1(q)|$ or $q^3$.
\end{Prf}

Finally, we consider the finite simple classical groups. 
Again, we do not check the converses of assertions~(i) or~(ii), remarking only that one could do so using the monograph~\cite{BG} of Burness and Giuidici, where the conjugacy classes of elements of prime order in these groups are classified.

\begin{La} \label{centralisersClassical}
Let $T$ be a finite simple classical group.
\begin{itemize}
\item[\textnormal{(i)}] If $T$ has type $\A_n^\varepsilon$, $\D_n$ or ${}^2\D_n$ with $n \ge 9$, or type $\B_n$ or $\C_n$ with $n \ge 5$, then there exists $x \in T \setminus \{1\}$ with $|C_T(x)| > |T|^{4/5}$.
\item[\textnormal{(ii)}] If $T$ has type $\A_n^\varepsilon$ with $n \ge 4$, type $\B_n$ or $\C_n$ with $n \ge 3$, or type $\D_n$ or ${}^2\D_n$ with $n \ge 4$, then there exists $x \in T \setminus \{1\}$ with $|C_T(x)| > |T|^{3/5}$.
\item[\textnormal{(iii)}] $|C_T(x)| < |T|^{2/5}$ for all $x \in T \setminus \{1\}$ if and only if $T \cong \PSL_2(q)$ with $q \neq 7$.
\end{itemize}
\end{La}

\begin{Prf}
Throughout the proof, we write $q = p^f$ with $p$ a prime and $f \ge 1$. 
First suppose that $T$ has type $\A_1$. 
That is, $T \cong \PSL_2(q)$, with $q \ge 4$. 
The smallest non-trivial conjugacy class of $T$ has size $q(q-1)$, $\tfrac{1}{2}(q^2-1)$ or $\tfrac{1}{2}q(q-1)$ according as whether $p=2$, $q \equiv 1 \pmod 4$ or $q \equiv 3 \pmod 4$. 
Hence, every non-trivial conjugacy class of $T$ has size greater than $|T|^{3/5}$ if and only if $q \neq 7$. 
Equivalently, $|C_T(x)| < |T|^{2/5}$ for all $x \in T \setminus \{1\}$ if and only if $q \neq 7$.

Now suppose that $T$ has type $\A_n^\varepsilon$ with $n \ge 2$. 
That is, $T \cong \PSL^\varepsilon_{n+1}(q)$ (where $\text{L}^+ := \text{L}$, $\text{L}^- := \text{U}$). 
Let $x \in G := \PGL^\varepsilon_{n+1}(q)$ be an element of order $p$ with one Jordan block of size $2$ and $n-1$ Jordan blocks of size $1$. 
That is, $a_1=n-1$, $a_2=1$ and $a_3=\cdots=a_p=0$ in the notation of \cite[Section~3.2.3]{BG}.
Then $x \in T$, and $x^T=x^G$ by \cite[Propositions~3.2.7 and~3.3.10]{BG}, so by \cite[Tables~B.3 and~B.4]{BG} we have $|C_T(x)| = |C_G(x)|/|G:T| = \tfrac{1}{d}|C_G(x)| = \tfrac{1}{d}q^{2n-1}|\GL^\varepsilon_{n-1}(q)|$, where $d := \gcd(n+1,q-\varepsilon)$. 
Therefore,
\begin{equation} \label{AnCentralisers1}
|C_T(x)| = \frac{1}{d}q^{n(n+1)/2}\prod_{i=1}^{n-1}(q^i-\varepsilon^i) \quad \text{and} \quad |T| = \frac{1}{d}q^{n(n+1)/2}\prod_{i=2}^{n+1}(q^i-\varepsilon^i).
\end{equation}
For $n \in \{2,3\}$ we must show that $|C_T(x)| > |T|^{2/5}$.
If $n=2$ then $d \le 3$, so $|C_T(x)| \ge \tfrac{1}{3}q^3(q-\varepsilon)$ while $|T| \le q^3(q^2-\varepsilon^2)(q^3-\varepsilon^3)$. 
This implies that $|C_T(x)| > |T|^{2/5}$ for all $q \ge 7$, and one may check directly that this inequality also holds for $q<7$. 
If $n=3$ then $d \le 4$, so $|C_T(x)| \ge \tfrac{1}{4}q^6(q-\varepsilon)(q^2-\varepsilon^2)$ while $|T| \le q^6(q^2-\varepsilon^2)(q^3-\varepsilon^3)(q^4-\varepsilon^4)$. 
This implies that $|C_T(x)| > |T|^{2/5}$ for all $q \ge 3$, and a direct calculation shows that this inequality also holds for $q=2$.
Now suppose that $4 \le n \le 8$. 
We must show that $|C_T(x)| > |T|^{3/5}$. 
Since $q \ge 2$, \eqref{AnCentralisers1} gives
\[
|C_T(x)| \ge \frac{1}{d} \frac{q^{n^2}}{2^{n-1}} \quad \text{and} \quad |T| \le \frac{1}{d} \Big(\frac{3}{2}\Big)^n q^{n^2+2n},
\]
and so it suffices to show that $q^{2n^2-6n} > d^2 2^{2n-5} 3^{3n}$. 
Indeed, since $d \le n+1$, it suffices to show that $q^{2n^2-6n} > (n+1)^2 2^{2n-5} 3^{3n}$. 
This inequality holds for all $q \ge 2$ if $n \in \{7,8\}$, for all $q \ge 3$ if $n=6$, for all $q \ge 4$ if $n=5$, and for all $q \ge 11$ if $n=4$. 
In the remaining cases, where $(n,q) = (6,2)$, $(5,2)$, $(5,3)$ or $(4,q)$ with $q<11$, one may check directly that $|C_T(x)| > |T|^{3/5}$.
It remains to show that $|C_T(x)| > |T|^{4/5}$ for all $q \ge 2$ when $n \ge 9$. 
If $q \ge 3$ then \eqref{AnCentralisers1} gives
\[
|C_T(x)| \ge \frac{1}{d} \Big(\frac{2}{3}\Big)^{n-1} q^{n^2} \quad \text{and} \quad |T| \le \frac{1}{d} \Big(\frac{4}{3}\Big)^n q^{n^2+2n},
\]
so it suffices to show that $q^{n^2-8n} > d \cdot 2^{3n+5} 3^{n-5}$. 
Indeed, since $d \le n+1$, we can just show that $q^{n^2-8n} > (n+1) 2^{3n+5} 3^{n-5}$. 
This inequality holds for all $q \ge 3$ if $n \ge 11$; if $n=10$, it holds for all $q \ge 5$, and if $n=9$, it holds for all $q \ge 29$. 
For $n=10$ with $2 \le q<5$ and $n=9$ with $2 \le q<29$, one may check directly that $|C_T(x)| > |T|^{4/5}$. 
Finally, we must check that $|C_T(x)| > |T|^{4/5}$ when $q=2$ and $n \ge 9$. 
Since $n \ge 9$, and since $\tfrac{255}{256} q^i \le q^i-\varepsilon \le \tfrac{257}{256} q^i$ for $i \ge 8$, \eqref{AnCentralisers1} gives
\[
|C_T(x)| \ge \frac{1}{d} \Big(\frac{255}{256}\Big)^{n-8} q^{n^2-28} \prod_{i=1}^7(q^i-1) \quad \text{and} \quad |T| \le \frac{1}{d} \Big(\frac{257}{256}\Big)^{n-8} q^{n^2+2n-27} \prod_{i=2}^7(q^i-1).
\]
Noting also that $d \le 3$, we see that it suffices to show that
\[
2^{n^2-8n-32} 255^{5n-40} \prod_{i=2}^7 (q^i-1) > 3 \cdot 256^{n-8} 257^{4n-32}.
\]
This inequality holds for all $n \ge 9$, and so the proof of the $\A_n^\varepsilon$ case is complete.

Next, suppose that $T$ has type $\C_n$, where $n \ge 2$. 
That is, $T \cong \PSp_{2n}(q)$. 
Write $G := \PGSp_{2n}(q)$, noting that $|G:T| = \gcd(2,q-1)$. 
If $p>2$, take $x \in G$ of order $p$ with one Jordan block of size $2$ and $2(n-1)$ Jordan blocks of size $1$. 
That is, $a_1 = 2(n-1)$, $a_2=1$ and $a_3=\cdots=a_p=0$ in the notation of \cite[Section~3.4.3]{BG}.
Then $x \in T$, and by \cite[Proposition~3.4.12]{BG}, $x^G$ splits into two $T$-conjugacy classes and hence $|C_T(x)| = 2|C_G(x)|/|G:T| = \tfrac{1}{2}q^{2n-1}|\Sp_{2(n-1)}(q)|$. 
If $p=2$ then $T=G$ and we take $x$ to be an involution of type $b_1$ as in \cite[Table~3.4.1]{BG}, so that $|C_T(x)| = q^{2n-1}|\Sp_{2(n-1)}(q)|$. 
Hence, for every $p$, we have 
\begin{equation} \label{CnCentralisers}
|C_T(x)| = \frac{1}{d} q^{n^2} \prod_{i=1}^{n-1}(q^{2i}-1) \quad \text{and} \quad |T| = \frac{1}{d} q^{n^2} \prod_{i=1}^n(q^{2i}-1), 
\end{equation}
where $d = \gcd(2,q-1) \le 2$. 
If $n=2$ then $|C_T(x)| \ge \tfrac{1}{2} q^4(q^2-1)$ and $|T| \le q^4(q^2-1)(q^4-1)$, and it follows that $|C_T(x)| > |T|^{2/5}$ for all $q \ge 2$. 
Similarly, for $n \in \{3,4\}$ one may check that $|C_T(x)| > |T|^{3/5}$ for all $q \ge 2$. 
Now suppose that $n \ge 5$. 
Since $q^{2i} \ge 4$, we have $q^{2i}-1 \ge \tfrac{3}{4}q^{2i}$ for all $i \ge 1$, and so $|C_T(x)| \ge \tfrac{1}{2} (\tfrac{3}{4})^{n-1} q^{2n^2-n}$, while $|T| < q^{2n^2+n}$. 
Hence, to show that $|C_T(x)| > |T|^{4/5}$, it suffices to show that $(\tfrac{3}{4})^{5n-5} q^{2n^2-9n} > 2$. 
This inequality holds for all $q \ge 2$ when $n \ge 6$, and for all $q \ge 4$ when $n=5$; for $(n,q) = (5,2)$ and $(5,3)$, one may check directly that $|C_T(x)| > |T|^{4/5}$. 

Now suppose that $T$ has type $\B_n$, where $n \ge 2$. 
That is, $T \cong \Omega_{2n+1}(q)$ with $q$ odd. 
For $q \equiv 1$ or $3 \pmod 4$, let $x \in G := \PGO_{2n+1}(q)$ be an involution of type $t_n$ or $t'_n$, respectively, in the sense of \cite[Sections~3.5.2.1 and~3.5.2.2]{BG}. 
Then $x \in T$ and $x^T = x^G$, so $|C_T(x)| = |C_G(x)|/|G:T| = \tfrac{1}{2}|C_G(x)| = |\SO_{2n}^\pm(q)|$ by \cite[Table~B.8]{BG}. 
Now, 
\begin{equation} \label{BnCentralisers}
|\SO_{2n}^\pm(q)| = q^{n^2-n}(q^n\mp 1)\prod_{i=1}^{n-1}(q^{2i}-1) > \frac{1}{2} q^{n^2} \prod_{i=1}^{n-1}(q^{2i}-1), 
\end{equation}
and the right-hand side above is the value of $|C_T(x)|$ that we obtained in the $\C_n$ case. 
Since $|\Omega_{2n+1}(q)| = |\PSp_{2n}(q)|$, we therefore reach the same conclusions as for type $\C_n$.

Now suppose that $T$ has type $\D_n^\varepsilon$, namely $T \cong \text{P}\Omega_{2n}^\varepsilon(q)$ with $n \ge 4$. 
Let $G := \text{Inndiag}(\text{P}\Omega_{2n}^\varepsilon(q))$, as defined on \cite[p.~56]{BG}. 
Assume first that $p>2$, noting that $|G:T|$ divides $4$. 
Take $x \in G$ of order $p$ with one Jordan block of size $2(n-2)$ and two Jordan blocks of size $2$. 
That is, $a_1=2(n-2)$, $a_2=2$ and $a_3=\cdots=a_p=0$ in the notation of \cite[Section~3.5.3]{BG}.
Then $x \in T$, and \cite[Propositions~3.5.14(i) and~(ii,b)]{BG} imply that $x^T = x^G$. 
Therefore, \cite[Table~B.12]{BG} gives $|C_T(x)| = |C_G(x)|/|G:T| \ge \tfrac{1}{4}|C_G(x)| = \tfrac{1}{8}q^{4n-7}|\text{O}_{2(n-2)}^{\varepsilon_1}(q)||\Sp_2(q)|$, where the value of $\varepsilon_1 = \pm$ depends on $n$ and $q$ as described there. 
Multiplying the inequality in \eqref{BnCentralisers} by $2$ to get a lower bound for $|\text{O}_{2(n-2)}^{\varepsilon_1}(q)|$, it follows that 
\begin{equation} \label{DnCentralisers1}
|C_T(x)| > \frac{1}{8} q^{n^2-2} (q^2-1) \prod_{i=1}^{n-3} (q^{2i}-1), \quad \text{while} \quad |T| = \frac{1}{d} q^{n(n-1)} (q^n-\varepsilon) \prod_{i=1}^{n-1}(q^{2i}-1),
\end{equation}
where $d = \gcd(4,q^n-\varepsilon)$. 
Since $q^{2i} \ge 9$ for all $i \ge 1$, and in particular $q^n \ge 3^4 = 81$, we have 
\begin{equation} \label{DnCentralisers2}
|C_T(x)| > \frac{1}{8} \Big(\frac{8}{9}\Big)^{n-2} q^{2n^2-5n+6} \quad \text{and} \quad |T| < \frac{82}{81} q^{2n^2-n}.
\end{equation}
For $4 \le n \le 8$ we need $|C_T(x)| > |T|^{3/5}$, so by \eqref{DnCentralisers2} it suffices to show that $(\tfrac{8}{9})^{5n-10} q^{4n^2-22n+30} > 8^5 (\tfrac{82}{81})^3$, which holds unless $(n,q) = (4,3)$ or $(4,5)$. 
For $(n,q) = (4,5)$, \eqref{DnCentralisers1} implies that $|C_T(x)| > |T|^{3/5}$; for $(n,q) = (4,3)$, a {\sf GAP} \cite{GAP} calculation shows that there exist elements $x \in T \setminus \{1\}$ for which this inequality holds. 
For $n \ge 9$ we claim that $|C_T(x)| > |T|^{4/5}$, and we now have $q^n \ge 3^9 = 19683$, so we can replace the $\tfrac{82}{81}$ in \eqref{DnCentralisers2} by $\tfrac{19684}{19683}$ to see that it suffices to show that $(\tfrac{8}{9})^{5n-10} q^{2n^2-21n+30} > 8^5 (\tfrac{19684}{19683})^4$. 
If $n \ge 10$ then this inequality holds for all $q \ge 3$, and if $n=9$ then it holds for $q \ge 127$. 
For $n=9$ and $q<127$, using the equality in \eqref{BnCentralisers} we obtain $|C_T(x)| \ge \tfrac{1}{4} q^{n^2-n} (q^{n-2}-1)(q^2-1)\prod_{i=1}^{n-3}(q^{2i}-1)$, which implies that $|C_T(x)| > |T|^{4/5}$ except when $q=3$ and $\varepsilon=+$. 
However, in this case we have $|G:T|=2$ (compare \cite[Figure~2.5.1 and Lemma~2.2.9]{BG}, noting that the discriminant of a hyperbolic quadratic form on $\F_q^{2n}$ with $(n,q)=(9,3)$ is $\boxtimes$, in the notation used there, because $n(q-1)/4=9$ is odd), so the $\tfrac{1}{4}$ in the above estimate for $|C_T(x)|$ may be replaced by $\tfrac{1}{2}$, and we again obtain $|C_T(x)| > |T|^{4/5}$.

Finally, suppose that $T \cong \text{P}\Omega_{2n}^\varepsilon(q)$ with $q$ even, noting that $T=G$ in this case. 
Take $x \in G$ to be an involution of type $a_2$ as in \cite[Table~3.5.1]{BG}. 
Then $|C_T(x)| = q^{4n-7}|\Omega_{2(n-2)}^\varepsilon(q)||\Sp_2(q)|$ and $\gcd(4,q^n-\varepsilon)=1$, so instead of \eqref{DnCentralisers1} we have
\begin{equation} \label{DnCentralisers3}
|C_T(x)| > \frac{1}{4} q^{n^2-2} (q^2-1) \prod_{i=1}^{n-3} (q^{2i}-1) \quad \text{and} \quad |T| = q^{n(n-1)} (q^n-\varepsilon) \prod_{i=1}^{n-1}(q^{2i}-1).
\end{equation}
(In the bound for $|C_T(x)|$ we drop a factor of $\tfrac{1}{4}$ because $|G:T|=1$ for $q$ even, but pick up a factor of $\tfrac{1}{2}$ because $\Omega_{2(n-2)}^\varepsilon(q)$ has index $2$ in $\SO_{2(n-2)}^\varepsilon(q)$.) 
For $4 \le n \le 8$ we need $|C_T(x)| > |T|^{3/5}$. 
Since $q^i \ge 2$ for all $i \ge 1$, and in particular $q^n \ge 16$, it suffices to show that $(\tfrac{3}{4})^{5n-10} q^{4n^2-22n+30} > 4^5 (\tfrac{17}{16})^3$. 
This inequality holds unless $(n,q) = (4,2)$ or $(4,4)$, in which cases a direct calculation shows that $|C_T(x)| > |T|^{3/5}$. 
For $n \ge 9$ we must show that $|C_T(x)| > |T|^{4/5}$. 
We now have $q^n \ge 512$, and so it suffices to show that $(\tfrac{3}{4})^{5n-10} q^{2n^2-21n+30} > 4^5 (\tfrac{513}{512})^4$. 
This inequality holds unless $(n,q)=(10,2)$ or $n=9$ and $q \le 2^8$. 
One may use \eqref{DnCentralisers3} to check that $|C_T(x)| > |T|^{4/5}$ in each of these cases except $(n,q)=(9,2)$, in which case the desired inequality may be verified by a direct calculation. 
\end{Prf}

\section{Quasiprimitive point actions of type SD or CD} \label{sec:SDCD}

We now apply Corollary~\ref{corForPrim} to permutation groups $N$ that arise as subgroups of certain types of primitive groups. 
In some cases, we are also able to treat quasiprimitive groups, namely those in which every non-trivial normal subgroup is transitive.
In this section, we consider the case where the group $N$ in Corollary~\ref{corForPrim} has a `diagonal' action. 
Specifically, we work under the following hypothesis.

\begin{hypothesis} \label{DiagonalHyp}
Let $T$ be a non-Abelian finite simple group, let $k \ge 2$, and write $H = T^k$. 
Let $\Omega = \{ (y_1,\ldots,y_{k-1},1) \mid y_1,\ldots,y_{k-1} \in T \} \le H$, and let $H$ act on $\Omega$ by 
\begin{equation} \label{diagonalAction}
(y_1,\ldots,y_{k-1},1)^{(x_1,\ldots,x_k)} = (x_k^{-1}y_1x_1,\ldots,x_k^{-1}y_{k-1}x_{k-1},1).
\end{equation}
Suppose that $N = H^r$ is a collineation group of a thick finite generalised quadrangle $\mathcal{Q} = (\mathcal{P},\mathcal{L},\inc)$ of order $(s,t)$, such that $\mathcal{P} = \Omega^r$ and $N$ has the product action on $\mathcal{P}$. 
\end{hypothesis}

This situation arises when $N$ is the socle of a primitive permutation group $G \le \Sym(\Omega)$ of type HS, HC, SD or CD. 
For type HS (respectively HC) we have $k=2$ and $r=1$ (respectively $r \ge 2$), $G$ has two minimal normal subgroups, each isomorphic to $T^r$, and the socle of $G$ is $T^r \times T^r$, which is isomorphic to $N$. 
For type SD (respectively CD) we have $k \ge 2$ and $r=1$ (respectively $r \ge 2$), and $G$ has a unique minimal normal subgroup, isomorphic to $T^{kr} \cong N$. 
Note that the notation $k$ and $r$ is consistent with that of Table~\ref{tab:Primitive2}. 
Of course, $G$ must (usually) satisfy certain other conditions \cite[Section~6]{CherylONS} in order to actually be primitive, but these conditions are not needed for the proof of Proposition~\ref{diagonalPrelim}. 
It suffices that there is a subgroup of the form $N$. 
In particular, we are also able to treat quasiprimitive groups~\cite[Section~12]{CherylONS}, because the (action of the) socle of $G$ is the same as in the respective primitive types. 
(Note that a quasiprimitive group of type HS or HC is necessarily primitive, but a quasiprimitive group of type SD or CD need not be primitive.)

Proposition~\ref{diagonalPrelim} shows, in particular, that the parameter $r$ in Hypothesis~\ref{DiagonalHyp} can be at most $3$.
As illustrated after Corollary~\ref{corForPrim}, the proof relies on the information about centraliser orders in non-Abelian finite simple groups given in Proposition~\ref{classes}. 
We also observe that when $r=3$, there always exists a solution $(s,t) = (|\Omega|-1,|\Omega|+1)$ of the equation $|\Omega|^3 = |\Omega|^r = |\mathcal{P}| = (s+1)(st+1)$, and this solution satisfies properties~(ii) and~(iii) of Lemma~\ref{lemmaBasics}. 
Hence, although we are unable to rule out the case $r=3$ completely, we verify that this `obvious' situation cannot occur. 

\begin{Pro} \label{diagonalPrelim}
If Hypothesis~\ref{DiagonalHyp} holds then $r \le 3$ and $|C_T(x)| < |T|^{1-r/5}$ for all $x \in T \setminus \{1\}$, and in particular $T$ must appear in Table~$\ref{tab:summary4}$. 
Moreover, if $r=3$ then $(s,t) \neq (|\Omega|-1,|\Omega|+1)$.
\end{Pro}

\begin{Prf}
Note first that $|\mathcal{P}| = |\Omega|^r = |T|^{(k-1)r}$.
In particular, the excluded case $(s,t) = (2,4)$ in Corollary~\ref{corForPrim} does not arise, because $|\mathcal{P}| \ge |T| \ge |\Alt_5| = 60 > (2+1)(2\cdot4+1)$. 
If we take $x := x_1=\cdots=x_k \neq 1$ in \eqref{diagonalAction}, then $(y_1,\ldots,y_{k-1},1) \in \Omega$ is fixed if and only if $y_1,\ldots,y_{k-1} \in C_T(x)$. 
That is, $(x,\ldots,x) \in H$ fixes precisely $|C_T(x)|^{k-1}$ elements of $\Omega$ (and, in particular, the action of $H$ on $\Omega$ is not semiregular). 
Corollary~\ref{corForPrim} therefore implies that $r \le 4$ and $|C_T(x)|^{k-1} < |\Omega|^{1-r/5} = |T|^{(k-1)(1-r/5)}$, namely $|C_T(x)| < |T|^{1-r/5}$, for all $x \in T \setminus \{1\}$. 
If $r=4$ then we have a contradiction because every non-Abelian finite simple group $T$ contains a non-identity element $x$ with $|C_T(x)| > |T|^{1/5}$. 
For example, it is well known that every non-Abelian finite simple group $T$ contains an involution $x$ with $|C_T(x)| > |T|^{1/3}$ (in fact, {\em every} involution has this property \cite[Proposition 2.4]{LiebeckShalevFixity}).  
Therefore, $r \le 3$. 
In particular, Proposition~\ref{classes} tells us that $T$ must be one of the groups appearing in Table~$\ref{tab:summary4}$.
To prove the final assertion, suppose towards a contradiction that $r=3$ and $(s,t) = (|\Omega|-1,|\Omega|+1)$. 
Take any $x \in T$ with $|C_T(x)| > |T|^{1/3}$. 
Then $((x,\ldots,x),(1,\ldots,1),(1,\ldots,1)) \in H^r = H^3 = N$ fixes $|C_T(x)|^{k-1}|T|^{2(k-1)} > |T|^{7(k-1)/3} = |\mathcal{P}|^{7/9}$ points of $\mathcal{Q}$, contradicting Proposition~\ref{t=s+2}.
\end{Prf}

The following immediate consequence of Proposition~\ref{diagonalPrelim} (and the preceding observations) implies the SD and CD cases of Theorem~\ref{thm:QP-GQs}.

\begin{Pro} \label{propDiagonalSummary}
Let $\mathcal{Q} = (\mathcal{P},\mathcal{L})$ be a thick finite generalised quadrangle admitting a collineation group $G$ that acts quasiprimitively of type SD or CD on $\mathcal{P}$. 
Then the conditions in Table~\ref{tab:Primitive2} hold. 
\end{Pro}

\section{Primitive point actions of type HS or HC} \label{sec:PA}

We now consider the case where $k=2$ in Hypothesis~\ref{DiagonalHyp} in more detail. 
As explained above, this case arises when $N$ is the socle of a primitive permutation group $G \le \Sym(\Omega)$ of type HS ($r=1$) or HC ($r \ge 2$). 
When $k=2$ it is natural to simplify the notation of Hypothesis~\ref{DiagonalHyp} by identifying the set $\Omega$ with $T^r$, so we first re-cast the hypothesis in this way and also establish some further notation.

\begin{hypothesis} \label{DiagonalHypk=2}
Let $T$ be a non-Abelian finite simple group and let $N = T^r \times T^r$ act on $T^r$ by 
\begin{equation} \label{diagonalActionk=2}
y^{(u_1,u_2)} = u_2^{-1}yu_1.
\end{equation}
Let $M = \{ (u,1) \mid u \in T^r \} \le N$, so that $M$ may be identified with $T^r$ acting regularly on itself by right multiplication.
Suppose that $N$ is a collineation group of a thick finite generalised quadrangle $\mathcal{Q} = (\mathcal{P},\mathcal{L},\inc)$ of order $(s,t)$ with $\mathcal{P} = T^r$. 
Let $\mathcal{P}_1 \subset \mathcal{P}$ denote the set of points collinear with but not equal to the identity element $1 \in T^r = \mathcal{P}$, and let $\mathcal{L}_1 \subset \mathcal{L}$ denote the set of lines incident with $1$.
Given a line $\ell \in \mathcal{L}$, let $\bar{\ell} \subset \mathcal{P}$ denote the set of points incident with $\ell$.
\end{hypothesis}

The following lemma may essentially be deduced from \cite[Lemma~10]{MR2287459} upon observing that the assumption $\gcd(s,t)>1$ imposed there is not necessary (as far as we can tell, and at least not in our more restrictive setting). 
We include a proof to make it clear that we do not need to make this assumption. 
Our notation differs from that of \cite[p.~654]{MR2287459} as follows: the point-regular group $G$ is our $M \cong T^r$, and the point $O$ is our point $1$, so that $\Delta$ is our $\mathcal{P}_1 \setminus \{1\}$.

\begin{La} \label{Yoshiara10}
Suppose that Hypothesis~\ref{DiagonalHypk=2} holds. 
Let $x \in \mathcal{P}_1 \setminus \{1\}$, and let $\ell_x$ be the unique line in $\mathcal{L}_1$ incident with $x$. 
Then, for every $i \in \{1,\ldots,|x|-1\}$, the conjugacy class $(x^i)^{T^r}$ is contained in $\mathcal{P}_1$. 
Moreover, the collineation $(x,1) \in M$ fixes $\ell_x$.
\end{La}

\begin{Prf}
Let us first establish some notation. 
Given $u \in T^r = \mathcal{P}$, write
\begin{align*}
\operatorname{fix}_\mathcal{P}(u) &= \{ P \in \mathcal{P} \mid P^{(u,1)} = P \}, \\
\operatorname{coll}_\mathcal{P}(u) &= \{ P \in \mathcal{P} \mid P^{(u,1)} \text{ is collinear with but not equal to } P \}, \\
\operatorname{fix}_\mathcal{L}(u) &= \{ \ell \in \mathcal{L} \mid \ell^{(u,1)} = \ell \}, \\
\operatorname{conc}_\mathcal{L}(u) &= \{ \ell \in \mathcal{L} \mid \ell^{(u,1)} \text{ is concurrent with but not equal to } \ell \}.
\end{align*}
Since the subgroup $M = \{ (u,1) \mid u \in T^r \}$ of $N$ acts regularly on $\mathcal{P}$, $\operatorname{fix}_\mathcal{P}(u)$ is empty. 
Moreover, $P \in \operatorname{coll}_\mathcal{P}(u)$ if and only if $P^{(P^{-1},1)}=1$ and $(Pu)^{(P^{-1},1)} = PuP^{-1}$ are collinear, which is if and only if $PuP^{-1} \in u^{T^r} \cap \mathcal{P}_1$. 
Since for $g,h \in T^r$ we have $gug^{-1}=huh^{-1}$ if and only if $g^{-1}h \in C_{T^r}(u)$, it follows that
\[
|\operatorname{coll}_\mathcal{P}(u)| = |u^{T^r} \cap \mathcal{P}_1||C_{T^r}(u)|,
\]
as in the proof of \cite[Lemma~3]{MR2287459}. 
Then (again, as in that proof) \cite[1.9.2]{FGQ} implies that
\begin{equation} \label{YoshiaraL3}
|\operatorname{coll}_\mathcal{P}(u)| = (s+1)|\operatorname{fix}_\mathcal{L}(u)| + |\operatorname{conc}_\mathcal{L}(u)| = |u^{T^r} \cap \mathcal{P}_1| |C_{T^r}(u)| \quad \text{(for every } u \in T^r\text{)}.
\end{equation}

Now, since $x \in \mathcal{P}_1$, we have $u^{-1}xu = x^{(u,u)} \in \mathcal{P}_1$ for every collineation of the form $(u,u) \in N$, because such collineations (are precisely those that) fix the point $1$. 
That is, every $T^r$-conjugate of $x$ is in $\mathcal{P}_1$. 
In other words, $x^{T^r} \cap \mathcal{P}_1 = x^{T^r}$, and so \eqref{YoshiaraL3} implies that
\begin{equation} \label{Benson}
|\operatorname{coll}_\mathcal{P}(x)| = (s+1)|\operatorname{fix}_\mathcal{L}(x)| + |\operatorname{conc}_\mathcal{L}(x)| = |x^{T^r}||C_{T^r}(x)| = |T^r| = |\mathcal{P}| = (s+1)(st+1).
\end{equation}
In particular, we have $\operatorname{coll}_\mathcal{P}(x) = \mathcal{P}$; that is, every point of $\mathcal{Q}$ is mapped to a collinear point under the collineation $(x,1) \in M$.
We now claim that $\operatorname{conc}_\mathcal{L}(x)$ is empty. 
If not, then some line $\ell$ is concurrent with its image under the collineation $(x,1)$. 
Let $P$ denote the unique point incident with both $\ell$ and $\ell^{(x,1)}$. 
Then $Px^{-1}$ is incident with $\ell$, being the image of $P$ under the collineation $(x,1)^{-1} = (x^{-1},1)$, and $Px^{-1} \neq P$ because $x \neq 1$ and $M$ acts regularly on $\mathcal{P}$. 
Since $\mathcal{Q}$ is thick, there exists a third point $P_3$ incident with $\ell$, distinct from $P$ and $Px^{-1}$. 
Since $\operatorname{coll}_\mathcal{P}(x) = \mathcal{P}$, the points $P_3^{(x,1)} = P_3x$ and $P_3$ are collinear. 
Moreover, $P_3x$ is collinear with $P$, because both of these points are incident with $\ell^{(x,1)}$. 
Therefore, $P_3x$ is collinear with two distinct points that are incident with $\ell$, namely $P_3$ and $P$, and so $P_3x$ is itself incident with $\ell$ because $\mathcal{Q}$ contains no triangles. 
This, however, means that $P_3x$ is incident with both $\ell$ and $\ell^{(x,1)}$, which forces $P_3x=P$ and hence $P_3=Px^{-1}$, a contradiction. 
Therefore, $|\operatorname{conc}_\mathcal{L}(x)|=0$ as claimed, and so \eqref{Benson} implies that
\begin{equation} \label{Benson2}
|\operatorname{fix}_\mathcal{L}(x)| = st+1.
\end{equation}
Next, we show that $(x^i)^{T^r} \subseteq \mathcal{P}_1$ for all $i \in \{1,\ldots,|x|-1\}$. 
For each such $i$, we certainly have $\operatorname{fix}_\mathcal{L}(x) \subseteq \operatorname{fix}_\mathcal{L}(x^i)$, because if the collineation $(x,1)$ fixes a line then so too does $(x,1)^i=(x^i,1)$. 
In particular, $|\operatorname{fix}_\mathcal{L}(x^i)| \ge |\operatorname{fix}_\mathcal{L}(x)| = st+1$, by \eqref{Benson2}. 
On the other hand, no two lines fixed by $(x^i,1)$ can be concurrent, because if they were, then the unique point incident with both lines would be fixed by $(x^i,1)$, a contradiction since $M$ acts regularly on $\mathcal{P}$. 
Hence, the total number of points that are incident with some line in $\operatorname{fix}_\mathcal{L}(x^i)$ is $(s+1)|\operatorname{fix}_\mathcal{L}(x^i)|$. 
As this number cannot exceed $|\mathcal{P}|=(s+1)(st+1)$, we must also have $|\operatorname{fix}_\mathcal{L}(x^i)| \le st+1$. 
Therefore, $|\operatorname{fix}_\mathcal{L}(x^i)| = st+1$. 
Now \eqref{YoshiaraL3} implies, on the one hand, that
\[
|\operatorname{coll}_\mathcal{P}(x^i)| = (s+1)|\operatorname{fix}_\mathcal{L}(x^i)| + |\operatorname{conc}_\mathcal{L}(x^i)| = |\mathcal{P}| + |\operatorname{conc}_\mathcal{L}(x^i)|.
\]
Since $|\operatorname{coll}_\mathcal{P}(x^i)| \le |\mathcal{P}|$, this implies that $|\operatorname{conc}_\mathcal{L}(x^i)| = 0$, and then in turn that $|\mathcal{P}| = |\operatorname{coll}_\mathcal{P}(x^i)|$. 
Appealing again to \eqref{YoshiaraL3}, we now deduce that $|(x^i)^{T^r} \cap \mathcal{P}_1| |C_{T^r}(x^i)| = |\mathcal{P}| = |T^r|$, which implies that $(x^i)^{T^r} \subseteq \mathcal{P}_1$ as required. 
The first assertion is therefore proved.

Finally, we must show that the collineation $(x,1)$ fixes the unique line $\ell_x \in \mathcal{L}_1$ incident with $x$. 
If $|x|=2$, then $(x,1)$ fixes $\ell_x$ because it fixes setwise the subset $\{1,x\}$ of points incident with $\ell_x$. 
That is, it maps $1$ to $1^{(x,1)}=1x=x$ and $x$ to $x^{(x,1)}=x^2=1$. 
Now suppose that $|x|>2$. 
Then the point $1^{(x,1)^2}=x^2 \neq 1$ is collinear with $x$ because $x$ is collinear with $1$. 
On the other hand, $(x^2)^{T_r} \subseteq \mathcal{P}_1$ by the first assertion, so in particular $x^2$ is collinear with $1$. 
Therefore, $x^2$ is collinear with two distinct points incident with $\ell_x$ (namely $1$ and $x$), and so is itself incident with $\ell_x$ because $\mathcal{Q}$ contains no triangles. 
Hence, $(x,1)$ fixes $\ell_x$ because it maps two points incident with $\ell_x$, namely $1$ and $x$, to another two points incident with $\ell_x$, namely $x$ and $x^2$.
\end{Prf}

Hypothesis~\ref{DiagonalHypk=2} imposes the following restrictions on the points and lines incident with the identity element of $T^r = \mathcal{P}$, and on the order $(s,t)$ of $\mathcal{Q}$. 

\begin{La} \label{lemmaLineSubgroup}
The following assertions hold under Hypothesis~\ref{DiagonalHypk=2}.
\begin{itemize}
\item[\textnormal{(i)}] $\mathcal{P}_1$ is a union of $T^r$-conjugacy classes.
\item[\textnormal{(ii)}] Every line $\ell \in \mathcal{L}_1$ has the property that $\bar{\ell}$ is a subgroup of $T^r$. 
Specifically, 
\[
\bar{\ell} = \{ u \in T^r \mid (u,1) \in M \text{ fixes } \ell \}.
\] 
\item[\textnormal{(iii)}] Every line $\ell \in \mathcal{L}_1$ is incident with an involution. 
\item[\textnormal{(iv)}] If some line in $\mathcal{L}_1$ is incident with representatives of every $T^r$-conjugacy class of involutions in $\mathcal{P}_1$, then $N$ acts transitively on the flags of $\mathcal{Q}$ and $r \ge 2$. 
\item[\textnormal{(v)}] $T^r$ has at least three conjugacy classes of involutions. 
\item[\textnormal{(vi)}] If $T^r$ has exactly three conjugacy classes of involutions, then either $\mathcal{P}_1$ contains exactly two of these classes, or $N$ acts transitively on the flags of $\mathcal{Q}$ and $r \ge 2$.
\item[\textnormal{(vii)}] $\gcd(s,t)=1$ and $t \ge s+1$.
\end{itemize}
\end{La}

\begin{Prf}
(i) This follows immediately from Lemma~\ref{Yoshiara10}.  

(ii) If $u \in \bar{\ell}$ then the collineation $(u,1) \in M$ fixes $\ell$ by Lemma~\ref{Yoshiara10}. 
Conversely, if $(u,1)$ fixes $\ell$ then, because $1 \in \mathcal{P}$ is incident with $\ell$, so too is $1^{(u,1)}=u$; that is, $u \in \bar{\ell}$. 

(iii) If $\ell \in \mathcal{L}_1$ is not incident with any involution, then $\bar{\ell}$, which is a subgroup of $T^r$ by (ii), must have odd order. 
That is, $s+1 = |\bar{\ell}|$ must be odd. 
However, $(s+1)(st+1) = |T|^r$ is even by the Feit--Thompson Theorem~\cite{MR0166261}, so $s$ must be odd and hence $s+1$ must be even, a contradiction. 

(iv) If $\ell \in \mathcal{L}_1$ is incident with representatives of every conjugacy class of involutions in $\mathcal{P}_1$, then $\ell$ can be mapped to any other line in $\mathcal{L}_1$ by some element of the stabiliser $N_1 = \{ (u,u) \mid u \in T^r \}$ in $N$ of the point $1 \in \mathcal{P} = T^r$. 
Since $N$ acts transitively on $\mathcal{P}$, this means that $N$ acts transitively on the flags of $\mathcal{Q}$. 
If $r=1$, this contradicts the main result of our earlier paper \cite[Theorem~1.1]{OurHSHC}, so $r \ge 2$.

(v) Suppose towards a contradiction that $T^r$ contains at most two conjugacy classes of involutions. 
Then $r=1$, because if $r \ge 2$ then any involution $x \in T$ gives rise to the three pairwise non-conjugate involutions $(x,1,\ldots,1)$, $(1,x,1,\ldots,1)$ and $(x,x,1,\ldots,1)$ in $T^r$. 
Hence, by (iii), $T$ must have exactly two conjugacy classes of involutions, say $x^T$ and $y^T$, and both must be contained in $\mathcal{P}_1$.
Without loss of generality, $x$ and $y$ commute, because at least one of them centralises a Sylow $2$-subgroup of $T$. 
Therefore, $xy$ is an involution, and so must be collinear with $1 \in \mathcal{P}$. 
Since $1$ is collinear with $x$, the images of $1$ and $x$ under the collineation $(y,1) \in M$ are collinear. 
That is, $1^{(y,1)}=y$ is collinear with $x^{(y,1)}=xy$. 
Similarly, $1$ and $y$ are collinear, and hence so too are $1^{(x,1)}=x$ and $y^{(x,1)}=yx=xy$. 
Since the involution $xy$ is also collinear with $1$ and $\mathcal{Q}$ contains no triangles, the points $1$, $x$, $y$ and $xy$ must be incident with a common line. 
In particular, $x$ and $y$ are incident with a common line in $\mathcal{L}_1$. 
Since $r=1$, this contradicts (iv).

(vi) Let $x$, $y$ and $z$ denote representatives of the three $T^r$-conjugacy classes of involutions. 
If $\mathcal{P}_1$ contains exactly one of these classes, then $N$ acts flag-transitively by (iv), and it follows from \cite[Theorem~1.1]{OurHSHC} that $r \ge 2$. 
Now suppose that $\mathcal{P}_1$ contains all three of $x^T$, $y^T$ and $z^T$. 
Without loss of generality, $x$ centralises a Sylow $2$-subgroup of $T^r$ and both $y$ and $z$ commute with $x$, so $xy=yx$ and $xz=zx$ are involutions. 
Arguing as in the proof of (iii), we deduce that $1$, $x$, $y$ and $xy$ are incident with a common line $\ell \in \mathcal{L}_1$. 
Replacing $y$ by $z$ in this argument, we see that $z$ is also incident with $\ell$, so (iv) again implies that $N$ acts flag-transitively (and it follows as above that $r \ge 2$).

(vii) If $\gcd(s,t)>1$ then \cite[Lemma~6(i)]{MR2287459} implies that {\em every} $T^r$-conjugacy class intersects $\mathcal{P}_1$. 
However, assertion~(i) then implies that $\mathcal{P}_1 = \mathcal{P}$, which is impossible. 
Therefore, $\gcd(s,t)=1$. 
In particular, to show that $t \ge s+1$ it suffices to show that $t \ge s$. 
The proof of this assertion is adapted from that of \cite[Corollary~2.3]{OurHSHC}. 
Choose two distinct lines $\ell_1,\ell_2 \in \mathcal{L}_1$, so that $\bar{\ell}_1$ and $\bar{\ell}_2$ are subgroups of $T^r$ by (ii). 
For brevity, we now abuse notation slightly and identify $\ell_1$ and $\ell_2$ with $\bar{\ell}_1$ and $\bar{\ell}_2$, respectively, dropping the `bar' notation.
Since $\ell_1$ is a subgroup of $T^r$ and right multiplication by any element of $T^r$ is a collineation of $\mathcal{Q}$ (identified with an element of $M$), we have in particular that every right coset $\ell_1 g_2$ of $\ell_1$ with $g_2 \in \ell_2$ corresponds precisely to the set of points incident with some line of $\mathcal{Q}$. 
Similarly, left multiplications (identified with elements of $\{1 \} \times T^r \le N$) are collineations, so every left coset $g_1 \ell_2$ of $\ell_2$ with $g_1 \in \ell_1$ is a line of $\mathcal{Q}$. 
Therefore, $\mathcal{L}' = \{ g_1\ell_2 \mid g_1 \in \ell_1\} \cup \{ \ell_1g_2 \mid g_2 \in \ell_2 \}$ is a subset of $\mathcal{L}$. 
Consider also the subset $\mathcal{P}'=\ell_1\ell_2$ of $\mathcal{P}=T^r$, and let $\inc'$ be the restriction of $\inc$ to $(\mathcal{P}' \times \mathcal{L}') \cup (\mathcal{L}' \times \mathcal{P}')$. 
If we can show that $\mathcal{Q}' = (\mathcal{P}',\mathcal{L}',\inc')$ is a generalised quadrangle of $\mathcal{Q}$ of order $(s,1)$, then \cite[2.2.2(i)]{FGQ} will imply that $t \ge s$. 
Let us first check that $\mathcal{Q}'$ satisfies the generalised quadrangle axiom. 
Let $\ell \in \mathcal{L}'$ and take $P \in \mathcal{P}'$ not incident with $\ell$. 
Then, since $\mathcal{Q}$ satisfies the generalised quadrangle axiom, there is a unique point $P_0 \in \mathcal{P}$ incident with $\ell$ and collinear with $P$. 
Since $\ell \subset \mathcal{P}'$, we have $P_0 \in \mathcal{P}'$, and so $\mathcal{Q}'$ also satisfies the generalised quadrangle axiom. 
It remains to check that $\mathcal{Q}'$ has order $(s,1)$.
Every line in $\mathcal{L}'$ is incident with $s+1$ points in $\mathcal{P}'$, being a coset of either $\ell_1$ or $\ell_2$, so it remains to show that every point in $\mathcal{P}'$ is incident with exactly two lines in $\mathcal{L}'$. 
Given $P=g_1g_2 \in \mathcal{P}'$, where $g_1 \in \ell_1$, $g_2 \in \ell_2$, each line $\ell \in \mathcal{L}'$ incident with $P$ is either of the form $h_1 \ell_2$ for some $h_1 \in \ell_1$ or $\ell_1 h_2$ for some $h_2 \in \ell_2$, and since $P \in \ell$, we must have $h_1=g_1$ or $h_2=g_2$, respectively. 
Therefore, $P$ is incident with exactly two lines in $\mathcal{L}'$, namely $g_1 \ell_2$ and $\ell_1 g_2$.
\end{Prf}

Proposition~\ref{diagonalPrelim} restricts the possibilities for the simple group $T$ in Hypothesis~\ref{DiagonalHypk=2} to those listed in Table~$\ref{tab:summary4}$. 
The following result shows that, furthermore, $T$ must be a Lie type group.

\begin{Pro} \label{diagonalPrelimk=2}
If Hypothesis~\ref{DiagonalHypk=2} holds then $T$ is a Lie type group.
\end{Pro}

\begin{Prf}
We have $|\mathcal{P}| = |T|^r$, and $r \in \{1,2,3\}$ by Proposition~\ref{diagonalPrelim}. 
For each of the alternating and sporadic simple groups $T$ in Table~$\ref{tab:summary4}$, we check computationally for solutions of $|T|^r = (s+1)(st+1)$ satisfying $s \ge 2$, $t \ge 2$ and properties~(ii) and~(iii) of Lemma~\ref{lemmaBasics} (see Remark~\ref{Rem:computations}).
If $r=3$ then the only such solutions have the form $(s,t) = (|T|-1,|T|+1) = (|\Omega|-1,|\Omega|+1)$, and this contradicts the final assertion of Proposition~\ref{diagonalPrelim}.
If $r \in \{1,2\}$ then the possibilities for $T$ and $(s,t)$ are as in Table~\ref{tab:k=2AltSpor}. 
By Lemma~\ref{lemmaLineSubgroup}(i), $\mathcal{P}_1$ is a union of $T^r$-conjugacy classes, and so we must be able to partition $|\mathcal{P}_1|=s(t+1)$ into a subset of the sizes of these classes (respecting multiplicities). 
When $r=1$ and $T \cong \Alt_7$ or $\Alt_8$, this is impossible: the non-trivial conjugacy class sizes not exceeding $s(t+1)$ are $70$, $105$ and $210$ in the first case, and $105$, $112$ and $210$ in the second (with each occurring exactly once). 
Similarly, if $r=2$ and $T\cong \text{J}_1$, one may check computationally that there is no partition of $s(t+1)$, where $(s,t)=(419,175141)$, into non-trivial $T^2$-conjugacy class sizes.
Hence, it remains to consider the cases where $r=2$ and $T \cong \Alt_6$ or $\text{M}_{11}$. 
Here we first determine computationally the possible partitions of $s(t+1)$ into non-trivial $T^2$-conjugacy class sizes to obtain a list of possible partitions $\mathcal{P}_1$ into $T^2$-conjugacy classes.
Now, because the point graph of $\mathcal{Q}$ is a strongly regular graph in which adjacent vertices have $\lambda:=s-1$ common neighbours and non-adjacent vertices have $\mu:=t+1$ common neighbours, $\mathcal{P}_1$ must be a partial difference set of $T^2$ with these parameters. 
That is, each non-identity element $y \in T^2$ must have exactly $\lambda$ representations of the form $y=z_i z_j^{-1}$ for $z_i,z_j \in \mathcal{P}_1$ if $y \in \mathcal{P}_1$, and exactly $\mu$ such representations if $y \not \in \mathcal{P}_1$. 
A computation verifies that this condition is violated for each of the partitions of $\mathcal{P}_1$ determined in the previous step. 
\end{Prf}

\begin{Rem} \label{Rem:computations}
\textnormal{
In the proof of Proposition~\ref{diagonalPrelimk=2}, and at several other points in Sections~\ref{sec:PA} and~\ref{ss:PA}, we need to check computationally whether certain positive integers $X$ can be equal to the number of points of a thick finite generalised quadrangle. 
That is, we check for integral solutions $(s,t)$ of the equation $(s+1)(st+1) = X$ subject to the constraints $s \ge 2$, $t \ge 2$, $s^{1/2} \le t \le s^2 \le t^4$ and $s+t \mid st(st+1)$ imposed by Lemma~\ref{lemmaBasics}. 
In Section~\ref{sec:PA}, $X$ has the form $|T|^m$ for some non-Abelian finite simple group $T$ and some $m \le 3$, and in Section~\ref{ss:PA} we instead have $X=Y^m$ with $m \le 4$ and $Y$ the index of a maximal subgroup of an almost simple group. 
The above inequalities imply that $s$ must lie between $X^{1/4}-1$ and $X^{5/2}$, so it suffices to consider every integer $s$ in this range and determine whether $t=((X-1)/s-1)/(s+1)$ is an integer and, if so, whether $s+t \mid st(st+1)$. 
We remark that we found it useful to also observe that $s$ must divide $X-1$, because it turned out that $X-1$ had only a very small number of divisors in many of the cases that we had to consider.
}
\end{Rem}

\begin{table}[!t]
\begin{small}
\begin{center}
\begin{tabular}{lllll}
\toprule
$r$ & $T$ & $s$ & $t$ & $s(t+1)$ \\
\midrule
$1$ & $\Alt_7$ & $11$ & $19$ & $220$ \\ 
$1$ & $\Alt_8$ & $19$ & $53$ & $1026$ \\ 
$2$ & $\Alt_6$ & $19$ & $341$ & $6498$ \\
$2$ & $\text{M}_{11}$ & $89$ & $7831$ & $697048$ \\
$2$ & $\text{J}_1$ & $419$ & $175141$ & $73384498$ \\
\bottomrule
\end{tabular}
\end{center}
\caption{\small Alternating and sporadic groups in the proof of Proposition~\ref{diagonalPrelimk=2}.}\label{tab:k=2AltSpor}
\end{small}
\end{table}

We now show that $r$ cannot equal $3$, and deduce some further restrictions on $T$ when $r \in \{1,2\}$.

\begin{Pro} \label{PropHC-rLeq3}
If Hypothesis~\ref{DiagonalHypk=2} holds then $r \le 2$ and $T$ is a Lie type group with the property that $|C_T(x)| < |T|^{1-2r/9}$ for all $x \in T \setminus \{1\}$.
\end{Pro}

\begin{Prf}
By Propositions~\ref{diagonalPrelim} and \ref{diagonalPrelimk=2}, $T$ is a Lie type group and $r \le 3$. 
We now show that $|C_T(x)| < |T|^{1-2r/9}$ for all $x \in T \setminus \{1\}$ and deduce from this that $r \neq 3$.
Suppose, towards a contradiction, that there exists $x \in T \setminus \{1\}$ with $|C_T(x)| \ge |T|^{1-2r/9}$. 
Define $w=(x,1,\ldots,1)\in T^r$ and let $\mathcal{Q}_\theta = (\mathcal{P}_\theta,\mathcal{L}_\theta)$ be the substructure of $\mathcal{Q}$ fixed by $\theta = (w,w) \in N_1$. 
Then $\mathcal{P}_\theta = C_T(x) \times T^{r-1}$, and hence 
\begin{equation} \label{HCr=3Bound7/9}
|\mathcal{P}_\theta| \ge |T|^{(1-2r/9)+(r-1)} = |T|^{7r/9} = |\mathcal{P}|^{7/9}.
\end{equation}
Proposition~\ref{Prop7/9} then says that either $s \ge t+3$, or $(s,t) \in \{ (2,4),(3,9) \}$. 
The first of these conditions contradicts Lemma~\ref{lemmaLineSubgroup}(vii); the second implies that $|T|^r = (s+1)(st+1) \in \{27,112\}$, which is impossible because $|T| \ge 60$.
Hence, every $x \in T \setminus \{1\}$ must satisfy $|C_T(x)| < |T|^{1-2r/9}$. 
For $r=3$ this says that $|C_T(x)| < |T|^{1/3}$ for all $x \in T \setminus \{1\}$, a contradiction because we can always find some $x$ with $|C_T(x)| > |T|^{1/3}$ (as noted in the proof of Proposition~\ref{diagonalPrelim}). 
Therefore, $r \le 2$.
\end{Prf}

Proposition~\ref{PropHC-rLeq3} allows us to further reduce the list of candidates for the simple group $T$ in Hypothesis~\ref{DiagonalHypk=2} in the remaining cases $r \in \{1,2\}$. 
Let us first consider $r=2$.

\begin{Cor} \label{diagonalPrelimk=2andr=3}
If Hypothesis~\ref{DiagonalHypk=2} holds with $r=2$ then $T$ is of Lie type $\A_1$, $\A_2^\varepsilon$, ${}^2\B_2$ or ${}^2\G_2$. 
In particular, $T$ has a unique conjugacy class of involutions. 
\end{Cor}

\begin{Prf}
The result is verified by straightforward calculations involving the bound on centraliser orders imposed by Proposition~\ref{PropHC-rLeq3}, but we include the details in Appendix~\ref{ss:app}.
\end{Prf}

We can now prove the HC case of Theorem~\ref{thm:QP-GQs}.

\begin{Thm} \label{ThmHC}
If $\mathcal{Q}$ is a thick finite generalised quadrangle with a collineation group $G$ that acts primitively on the point set $\mathcal{P}$ of $\mathcal{Q}$, then the action of $G$ on $\mathcal{P}$ does not have O'Nan--Scott type HC.
\end{Thm}

\begin{Prf}
As explained above, the socle of $G$ is a group $N=T^r \times T^r$ as in Hypothesis~\ref{DiagonalHypk=2}, for some $r \ge 2$. 
However, Corollary~\ref{diagonalPrelimk=2andr=3} tells us that $r=2$ and that $T$ has a unique conjugacy class of involutions. 
In particular, $T^r = T^2$ has exactly three conjugacy classes of involutions, with representatives $(x,1)$, $(1,y)$ and $(x,y)$, where $x$ and $y$ are involutions in $T$. 
Now, \cite[Theorem~1.1]{OurHSHC} says that $G$ cannot act transitively on the flags of $\mathcal{Q}$, so in particular $N$ cannot act transitively on the flags of $\mathcal{Q}$. 
Lemma~\ref{lemmaLineSubgroup}(vi) therefore implies that $\mathcal{P}_1 = T^2$ must contain exactly two $T^2$-conjugacy classes. 
Hence, without loss of generality, $\mathcal{P}_1$ contains the class $(x,1)^{T^2}$. 
Since $G$ acts primitively on $\mathcal{P}$, it induces a subgroup of $\Aut(T^2) = \Aut(T) \wr \Sym_2$ that swaps the two simple direct factors of $T^2$. 
Therefore, $\mathcal{P}_1$ also contains the class $(1,y)^{T^2}$, and so does not contain the class $(x,y)^{T^2}$. 
In particular, no line $\ell \in \mathcal{L}_1$ can be incident with both a conjugate of $(x,1)$ and a conjugate of $(1,y)$, because by Lemma~\ref{lemmaLineSubgroup}(ii), $\ell$ would then also be incident with the product of these elements, a conjugate of $(x,y)$. 
Hence, $\mathcal{L}_1$ is partitioned into two sets of lines: those incident with conjugates of $(x,1)$, and those incident with conjugates of $(1,y)$. 
Since $G_1$ swaps these sets, $G$ acts flag-transitively, in contradiction with \cite[Theorem~1.1]{OurHSHC}.
\end{Prf}

For $r=1$ we are left with the following list of candidates for $T$.

\begin{Cor} \label{HSleftoverT}
Suppose that Hypothesis~\ref{DiagonalHypk=2} holds with $r=1$. 
Then $T$ is of Lie type $\A_1$, $\A_n^\varepsilon$ with $2 \le n \le 6$, $\B_2$, $\C_2$, $\C_3$, $\D_n^\varepsilon$ with $4 \le n \le 6$, or exceptional Lie type other than $\E_8$.
\end{Cor}

\begin{Prf}
By Propositions~\ref{diagonalPrelim} and~\ref{diagonalPrelimk=2}, $T$ is one of the Lie type groups in the first column of Table~\ref{tab:summary4}. 
By arguing as in the proof of Proposition~\ref{PropHC-rLeq3} but applying Proposition~\ref{Prop.752} instead of Proposition~\ref{Prop7/9} in the first paragraph, we conclude that one of the following conditions must also hold: 
\begin{itemize}
\item[(i)] every non-identity element $x \in T$ satisfies $|C_T(x)| < |T|^{94/125}$, or 
\item[(ii)] $s \le 2.9701 \times 10^{15}$.
\end{itemize}
By choosing appropriate elements $x \in T$ as in the proofs in Section~\ref{sec:groups}, we are able to use this to deduce that $T$ does not have type $\A_7^\varepsilon$, $\A_8^\varepsilon$, $\B_4$, $\C_4$, $\D_7^\varepsilon$, $\D_8^\varepsilon$ or $\E_8$. 
We rule out $\E_8$ here as an example, and include details of the remaining cases in Appendix~\ref{ss:app}.
If $T \cong \E_8(q)$ then $(s+1)^4 > |\mathcal{P}| = |T| \ge |\E_8(2)| \approx 3.378 \times 10^{74}$ and hence $s > |\E_8(2)|^{1/4}-1 \approx 4.287 \times 10^{18}$, contradicting (ii), so (i) must hold. 
However, as noted in the proof of Lemma~\ref{centralisersExceptional}, there exists $x \in T \cong \E_8(q)$ with $|C_T(x)| = q^{57}|\E_7(q)| \sim q^{190}$, while $|T|^{94/125} \sim (q^{248})^{94/125} < q^{187}$. 
Indeed, one may check that $|C_T(x)| > |T|^{94/125}$ for all $q \ge 2$.
\end{Prf}

Finally, we use Lemma~\ref{lemmaLineSubgroup} to reduce the list of candidates for $T$ in Corollary~\ref{HSleftoverT} to those given in the first row of Table~\ref{tab:Primitive2}, thereby proving the HS case of Theorem~\ref{thm:QP-GQs}. 

\begin{Pro} \label{HS-2B2}
Let $\mathcal{Q} = (\mathcal{P},\mathcal{L})$ be a thick finite generalised quadrangle admitting a collineation group $G$ that acts primitively of type HS on $\mathcal{P}$, with socle $T \times T$ for some non-Abelian finite simple group $T$. 
Then $T$ has Lie type $\A_5^\varepsilon$, $\A_6^\varepsilon$, $\B_3$, $\C_2$, $\C_3$, $\D_4^\varepsilon$, $\D_5^\varepsilon$, $\D_6^\varepsilon$, $\E_6^\varepsilon$, $\E_7$ or $\FF_4$. 
\end{Pro}

\begin{Prf}
We are assuming that Hypothesis~\ref{DiagonalHypk=2} holds with $r=1$, so $T$ must be one of the groups listed in Corollary~\ref{HSleftoverT}. 
It remains to show that, further, $T$ cannot have Lie type $\A_1$, $\A_2^\varepsilon$, $\A_3^\varepsilon$, $\A_4^\varepsilon$, ${}^2\B_2$, ${}^2\G_2$, ${}^2\FF_4$, $\G_2$ or ${}^3\D_4$. 
This follows from Lemma~\ref{lemmaLineSubgroup}(v), because in each of these cases $T$ has at most two conjugacy classes of involutions. 
(This may be verified using, for example, \cite[Table~4.5.1]{GLS3} for odd characteristic and \cite{LiebeckSeitz} for even characteristic.)
\end{Prf}

\begin{Rem} \label{HSremark}
\textnormal{
Proposition~\ref{HS-2B2} begs the obvious question of whether we can rule out the last remaining candidates for $T$ listed there. 
We are confident that we {\em will} eventually be able to do this, but it seems that it will require even more new ideas and a detailed case-by-case analysis. 
Of course, some of the remaining groups can be ruled out in certain cases using Lemma~\ref{lemmaLineSubgroup}(v); in particular, if $T$ has Lie type $\C_2$, $\C_3$, $\FF_4$ or $\E_6^\varepsilon$ in characteristic $p$, then we must have $p=2$, because in odd characteristic these groups have only two conjugacy classes of involutions. 
When $T$ has exactly three conjugacy classes of involutions, we can begin by applying Lemma~\ref{lemmaLineSubgroup}(v) (because we know from \cite{OurHSHC} that $N$ cannot act transitively on the flags of $\mathcal{Q}$), and then the arguments in the proof of Lemma~\ref{lemmaLineSubgroup}(iv--vi) can be extended to deduce some restrictions on {\em which} involutions can appear in $\mathcal{P}_1$. 
However, even with this extra information, we have thus far been unable to completely rule out any of the remaining candidates for $T$. 
These kinds of arguments become more difficult when $T$ has more than three conjugacy classes of involutions, and in any case, it seems that it will be necessary to treat each group individually, and to use the structure of its involution centralisers in some detail. 
Although not an ideal state of affairs, we therefore leave the remaining cases for a future project.
}
\end{Rem}

\section{Primitive point actions of type PA} \label{ss:PA}
We now apply Corollary~\ref{corForPrim} to the case where $N$ is the socle of a primitive permutation group $G$ of O'Nan--Scott type PA. 
The notation of Hypothesis~\ref{PAhyp} coincides with that of Table~\ref{tab:Primitive2}. 

\begin{hypothesis} \label{PAhyp}
Let $\mathcal{Q} = (\mathcal{P},\mathcal{L})$ be a thick finite generalised quadrangle of order $(s,t)$ admitting a collineation group $G$ that acts primitively of type PA on $\mathcal{P}$, writing $T^r \le G \le H \wr \Sym_r$ for some almost simple primitive group $H \le \Sym(\Omega)$ with socle $T$, where $r \ge 2$. 
\end{hypothesis}

\begin{Pro} \label{propPAcor}
If Hypothesis~\ref{PAhyp} holds then $2 \le r \le 4$ and every non-identity element of $H$ fixes less than $|\Omega|^{1-r/5}$ points of $\Omega$.
\end{Pro}

\begin{Prf}
The socle of $G$ is $N = T^r$ and the action of $H$ on $\Omega$ is not semiregular, so the result follows immediately from Corollary~\ref{corForPrim}.
\end{Prf}

To say more than this, we would like to have generic lower bounds for the {\em fixity} $f(H)$, namely the maximum number of fixed points of a non-identity element, of an almost simple primitive group $H \le \Sym(\Omega)$. 
This problem was investigated in a recent paper of Liebeck and Shalev~\cite{LiebeckShalevFixity}, who proved that $f(H) \ge |\Omega|^{1/6}$ except in a short list of exceptions. 
This lower bound is not quite large enough to force further restrictions on $r$ in Proposition~\ref{propPAcor}, because to rule out $r=4$ (as we did for types HC and CD) we would need $f(H)$ to be at least $|\Omega|^{1/5}$. 
However, Liebeck and Shalev remark (after \cite[Theorem~4]{LiebeckShalevFixity}) that their $|\Omega|^{1/6}$ bound could potentially be improved generically to around $|\Omega|^{1/3}$, which would be sufficient for this purpose. 
Work in this direction is currently being undertaken by Elisa Covato at the University of Bristol as part of her PhD research \cite{ElisaPhD}, with the aim of classifying the almost simple primitive permutation groups $H \le \Sym(\Omega)$ containing an involution that fixes at least $|\Omega|^{4/9}$ points. 
As of this writing, the alternating and sporadic cases have been completed, and so we are able to apply these results to sharpen Proposition~\ref{propPAcor} as follows.

\begin{Pro} \label{FixityAltn}
Suppose that Hypothesis~\ref{PAhyp} holds with $r>2$ and $T$ an alternating group or a sporadic simple group, and let $S \le H$ denote the point stabiliser in the action of $H$ on $\Omega$. 
Then $r \in \{3,4\}$, $H=T \cong \Alt_p$ with $p$ a prime congruent to $3$ modulo $4$, and $S \cap T = p.\frac{p-1}{2}$. 
\end{Pro}

\begin{Prf} 
Since $r>2$, Proposition~\ref{propPAcor} tells us that $r \in \{3,4\}$, and that the fixity $f(H)$ of $H$ must be at most $|\Omega|^{1-r/5}$. 
If $f(H) \le |\Omega|^{1-3/5} = |\Omega|^{2/5}$ then Covato's results \cite{ElisaPhD} imply that either (i) $T \cong \Alt_p$ with $p \equiv 3 \pmod 4$ a prime and $S \cap T = p.\frac{p-1}{2}$, or (ii) $H$ and $S$ are in Table~\ref{tab:PAElisa}. 

In case~(i) we can at least deal with the situation where $H = \Sym_p$. 
Indeed, by the argument in \cite[Section~6]{LiebeckShalevFixity}, there is an involution $u \in S = p.(p-1)$ fixing $2^{(p-3)/2} (\frac{p-3}{2})!$ elements of $\Omega$, which is greater than $|\Omega|^{2/5} = (2(p-2)!)^{2/5}$ provided that $p > 7$. 
If $p=7$ then we observe that $u$ still fixes more than $|\Omega|^{1/3}$ elements. 
This rules out $r=4$, because then $1/3 > 1-r/5 = 1/5$. 
For $r=3$ we apply Proposition~\ref{t=s+2}. 
We have $|\Omega| = 2\cdot 5! = 120$ and hence $|\mathcal{P}| = |\Omega|^3 = 120^3$, and the only solution of $120^3 = (s+1)(st+1)$ satisfying $s \ge 2$, $t \ge 2$ and properties~(ii) and~(iii) of Lemma~\ref{lemmaBasics} is $(s,t)=(119,121)$, so Proposition~\ref{t=s+2} implies that every non-identity collineation of $\mathcal{Q}$ fixes at most $|\mathcal{P}|^{7/9}$ points. 
However, the collineation $(u,1,1) \in G$ fixes more than $|\Omega|^{1/3}|\Omega|^2 = |\mathcal{P}|^{7/9}$ points, a contradiction. 
Therefore, if we are in case~(i) then we must have $H = \Alt_p$, as per the assertion.

Now suppose that $H$ and $S$ are in Table~\ref{tab:PAElisa}. 
First consider the six cases on the left-hand side of the table. 
In each of these cases, $f(H)$ is at least $|\Omega|^{1/5}$, so $r=4$ is ruled out. 
For $r=3$, we apply Proposition~\ref{t=s+2} as above. 
Since $f(H) > |\Omega|^{1/5}$, we have in particular $f(H) > |\Omega|^{1/6}$. 
Choose $u \in H$ fixing at least $|\Omega|^{1/6}$ elements of $\Omega$, and consider the collineation $(u,1,1) \in G$, which fixes at least $|\Omega|^{1/6+2} = |\Omega|^{13/6} = |\mathcal{P}|^{13/18}$ points of $\mathcal{Q}$. 
Since the only solutions of $|\Omega|^3 = (s+1)(st+1)$ satisfying $s \ge 2$, $t \ge 2$ and properties~(ii) and~(iii) of Lemma~\ref{lemmaBasics} are those with $t=s+2$, Proposition~\ref{t=s+2} provides a contradiction. 
Now consider the five cases on the right-hand side of Table~\ref{tab:PAElisa}. 
The actions of $\text{J}_3.2$, $\text{O'N}.2$ and $\text{Th}$ all have fixity greater than $|\Omega|^{1/6}$ \cite[Lemma~5.3]{LiebeckShalevFixity}, so these are ruled out for both $r=3$ and $r=4$ by the same arguments as above. 
Now consider the action of $\text{M}_{23}$. 
Here $|\Omega| = 40320$, and for $r=4$ there are no solutions of $|\Omega|^r = (s+1)(st+1)$ satisfying $s \ge 2$, $t \ge 2$ and properties~(ii) and~(iii) of Lemma~\ref{lemmaBasics}. 
If $r=3$, the only solution is $(s,t)=(40319,40321)$. 
By \cite[Lemma~5.3]{LiebeckShalevFixity}, we have $f(H)=5$, realised by an element $u$ of order $11$, and so we can construct a collineation $\theta=(u,1,1) \in G$ of $\mathcal{Q}$ fixing $5|\Omega|^2 = 8128512000$ points.
However, $s=t-2<t+3$, so the final assertion of Lemma~\ref{substructureLemma} implies that $|\mathcal{P}_\theta| \le  (s+1)(s+3) = 1988752683 < 8128512000$, and we have a contradiction.
Finally, consider the given action of $\text{B}$, for which $|\Omega|=3843461129719173164826624000000$. 
For $r=4$ there is no admissible solution of $|\Omega|^r = (s+1)(st+1)$. 
For $r=3$ the only admissible solution is $(s,t)=(|\Omega|-1,|\Omega|+1)$, and so the final assertion of Lemma~\ref{substructureLemma} implies that any non-identity collineation of $\mathcal{Q}$ fixes at most $(|\Omega|+1)(|\Omega|+3)$ points. 
However, \cite[Lemma~5.3]{LiebeckShalevFixity} tells us that $f(H) = 22$, so we can construct a collineation with $22|\Omega|^2$ points to yield a contradiction.
\end{Prf}

\begin{Rem} \label{RemPA}
\textnormal{
Further improvements to Proposition~\ref{propPAcor} will be made in a future project. 
In the first instance, we hope to use Covato's results \cite{ElisaPhD} on fixities of Lie type groups (once available), to complete our treatment of the cases $r=3$ and $r=4$. 
We also note that it is straightforward to handle the case $r=2$ with $T$ a sporadic simple group, and likewise the almost simple case with sporadic socle, computationally along the lines of \cite[Section~6]{Bamberg:2012yf} (but assuming only point-primitivity and not line-primitivity). 
However, we omit these computations from the present paper for brevity.
}
\end{Rem}

\begin{table}[!t]
\begin{small}
\begin{center}
\begin{tabular}{lllll}
\toprule
$H$ & $S$ && $H$ & $S$ \\
\midrule
$\Alt_9$ & $3^2:\SL_2(3)$ && $\text{J}_3.2$ & $19.9$ \\
$\text{J}_1$ & $2^3.7.3$ && $\text{O'N}.2$ & $31.15$ \\
$\text{J}_1$ & $7:6$ && $\text{M}_{23}$ & $23.11$ \\
$\text{He}$ & $7^2:2.\PSL_2(7)$ && $\text{Th}$ & $31.15$ \\
$\text{He}.2$ & $7^2:2.\PSL_2(7).2$ && $\text{B}$ & $47.23$ \\
$\text{Th}$ & $2^5.\PSL_5(2)$ &&& \\
\bottomrule
\end{tabular}
\end{center}
\caption{\small Actions with small fixity in Proposition~\ref{FixityAltn}.}\label{tab:PAElisa}
\end{small}
\end{table}

\section{Proof of Theorem~\ref{thm:QP-GQs}} \label{ss:mainProof}

Let us now summarise the proof of Theorem~\ref{thm:QP-GQs}. 
If the primitive action of $G$ on $\mathcal{P}$ has O'Nan--Scott type AS or TW, then the conditions stated in Table~\ref{tab:Primitive2} follow immediately from Theorem~\ref{subGQLemma-general}.
Types HS, HC and PA are treated in Proposition~\ref{HS-2B2}, Theorem~\ref{ThmHC} and Proposition~\ref{FixityAltn}, respectively. 
Types SD and CD are treated together in Proposition~\ref{propDiagonalSummary}.

\section{Discussion and open problems} \label{ss:discussion}

We feel that the results presented in this paper represent a substantial amount of progress towards the classification of point-primitive generalised quadrangles, but there is evidently still a good deal of work to do. 
We conclude the paper with a brief discussion, and outline some open problems that could be investigated independently and then potentially applied to our classification program. 

As discussed in Remark~\ref{HSremark}, we are confident that we will eventually be able to finish the HS case, and it is at least somewhat clear how this might be done.  
The SD and CD cases would appear to be more difficult, however. 
The arguments used in Section~\ref{sec:PA} do not work in these cases, because the proof of Lemma~\ref{Yoshiara10} (and therefore Lemma~\ref{lemmaLineSubgroup}) relies in a crucial way on having $k=2$ in Hypothesis~\ref{DiagonalHyp}, so that conjugation by an element of the underlying point-regular group $M$ is a collineation. 
We have thus far been unable to find a way to work around this difficulty in any sort of generality. 
On the other hand, a primitive (respectively quasiprimitive) group of type SD must induce a primitive (respectively transitive) permutation group on the set of simple direct factors of its socle $T^k$, and it seems that it should be possible to use this extra structure to say more about the SD and CD types, at least in the primitive case (especially since we have already reduced the list of candidates for $T$ to those in Table~\ref{tab:Primitive2}). 
Although we have made some preliminary investigations along these lines, we do not yet know how to finish the SD and CD cases, and so we leave this task for a future project.

\subsection{Point-regular collineation groups, and a number-theoretic problem}

There is, of course, a potential --- but apparently extremely challenging --- way to deal with all of the types HS, SD and CD, and also with type TW, in one fell swoop. 
In each of these cases, the full collineation group must have a point-regular subgroup of the form $T^m$, for some $m$, with $T$ a non-Abelian finite simple group. 
Hence, it would certainly be sufficient to show that such a group cannot act regularly on the points of a generalised quadrangle. 
However, this would appear to be a very difficult problem in light of the (limited) existing literature. 
Yoshiara~\cite{MR2287459} managed to show that a generalised quadrangle of order $(s,t)$ with $s=t^2$ cannot admit a point-regular group, while Ghinelli~\cite{MR1153980} considered the case where $s$ is even and $t=s$, showing that such a group must have trivial centre and cannot be a Frobenius group. 
Beyond this, not much else seems to be known in the way of restrictions on groups that can act regularly on the points of a generalised quadrangle (though certainly many of the {\em known} generalised quadrangles admit point-regular groups \cite{BG-regular}, and the Abelian case is understood \cite{MR2201385}). 
Although Yoshiara~\cite{MR2287459} has an extensive suite of lemmas that one might attempt to use to investigate (in particular) the possibility that a group of the form $T^n$ acts point-regularly on a generalised quadrangle $\mathcal{Q}$, the bulk of these lemmas assume that the order $(s,t)$ of $\mathcal{Q}$ satisfies $\gcd(s,t) \neq 1$. 
Although this condition holds under Yoshiara's intended assumption that $s=t^2$, it seems to be difficult to guarantee in general. 
Indeed, according to Lemma~\ref{lemmaLineSubgroup}(ii) (and perhaps not surprisingly), it must {\em fail} in our HS case. 
On the other hand, one might seek a contradiction by examining the arithmetic nature of the equation $|T|^m = (s+1)(st+1)$ subject to the constraints $s \ge 2$, $t \ge 2$, $s^{1/2} \le t \le s^2 \le t^4$ and $s+t \mid st(st+1)$ imposed by Lemma~\ref{lemmaBasics}, and asking when it can be guaranteed that a solution must satisfy $\gcd(s,t) \neq 1$. 
More generally, one might simply ask whether this equation can have any such solutions at all. 
This motivates the following problem.

\begin{Prob} \label{problem:sols}
Determine for which non-Abelian finite simple groups $T$, and which positive integers $m$, there exist integral solutions $(s,t)$ of the equation 
\begin{equation} \label{eq:sols}
|T|^m = (s+1)(st+1) \quad \text{with} \quad s \ge 2,\; t \ge 2,\; s^{1/2} \le t \le s^2 \le t^4 \text{ and } s+t \mid st(st+1).
\end{equation}
Failing this, determine when such a solution must satisfy $\gcd(s,t) \neq 1$.
\end{Prob}

As noted before Proposition~\ref{diagonalPrelim}, there is always an `obvious' solution of \eqref{eq:sols} when $m$ is divisible by $3$, namely $(s,t) = (|T|^{m/3}-1,|T|^{m/3}+1)$, and $\gcd(s,t)=1$ in this case because $|T|$ is even. 
It would be useful even to know whether this is the {\em unique} solution in this particular situation. 
We do know that \eqref{eq:sols} has solutions for certain $T$ when $m=1$ or $2$, as demonstrated by Table~\ref{tab:k=2AltSpor}, but we do not recall encountering any solutions apart from the aforementioned `obvious' ones when $m \ge 3$. 
Moreover, it is straightforward to run numerical computations that suggest that certain combinations of families of $T$ and values of $m$ will never yield a solution of \eqref{eq:sols}. 
For example, if $T \cong \PSL_2(q)$ and $m=1$ then there is no solution if $q < 10^6$, but we do not see how to go about proving that there is no solution for {\em any} $q$. 

One might also ask about gearing Problem~\ref{problem:sols} towards the PA and AS cases, by seeking solutions of \eqref{eq:sols} with $|T|$ replaced by $|H:S|$ for $H$ an almost simple group with socle $T$ and $S$ a maximal subgroup of $H$ (compare Hypothesis~\ref{PAhyp}, which reduces to the AS case if $r$ is taken to be $1$). 
However, solutions of \eqref{eq:sols} seem to be rather more common in this setting, and so other methods are needed to rule out cases where solutions arise. 
For example, if we take $H=T=\text{McL}$ (the McLaughlin sporadic simple group) then there are five (classes of) maximal subgroups $S$ of $H$ for which $|H:S|^2 = (s+1)(st+1)$ has an `admissible' solution: four maximal subgroups of order $40320$, which yield $(s,t)=(296,5644)$, and the maximal subgroup $\PSU_4(3)$, for which $(s,t)=(24,126)$. 

\subsection{Fixities of primitive groups of type TW}

We conclude with a brief discussion of the TW case. 
Let $N = T_1 \times \cdots \times T_r$, where $T_1 \cong \cdots \cong T_r \cong T$ for some non-Abelian finite simple group $T$. 
A primitive permutation group $G \le \Sym(\Omega)$ of type TW is a semidirect product $N \rtimes P$ with socle $N$ acting regularly by right mutlitplication, and $P \le \Sym_r$ acting by conjugation in such a way that $T_1,\ldots,T_r$ are permuted transitively. 
Certain other rather complicated conditions must also hold \cite{Baddeley}, and in particular $T$ must be a section of $P$.
If we intend to apply Theorem~\ref{thm:QP-GQs} to classify the generalised quadrangles with a point-primitive collineation group of TW type, then we will need `good' lower bounds for fixities of primitive TW-type groups. 
Liebeck and Shalev \cite[Section~4]{LiebeckShalevFixity} show that, for every $T$ and $r$, the fixity of $G$ is at least $|T|^{r/3}$. 
Although this is very far away from the $4/5$ exponent bound imposed by Theorem~\ref{thm:QP-GQs}, we would be interested to know under what conditions it could be improved to something `close' to $4/5$, so that we could at least rule out some of the subgeometries listed in Lemma~\ref{substructureLemma} and then perhaps use the underlying point-regular group to say more. 
In \cite[Section~4]{LiebeckShalevFixity}, Liebeck and Shalev consider an involution $x \in P$ (which must exist because $T$ is a section of $P$ and $|T|$ is even) and observe that $x$ induces a permutation of $\{T_1,\ldots,T_r\}$ that fixes at least $|T|^{ca+b}$ elements of $\Omega \equiv T^r$ , where the induced permutation has cycle structure $(1^a,2^b)$ and {\em every} involution $g \in \Aut(T)$ satisfies $|C_T(g)| \ge |T|^c$. 
By \cite[Proposition~2.4]{LiebeckShalevFixity}, we can take $c=1/3$ independently of $T$, and so because $a/3 + b \ge (r-2)/3+2/3=r/3$, it follows that $x$ fixes at least $|T|^{r/3}$ elements. 
Now, $c$ can certainly be taken larger than $1/3$ in at least some non-Abelian finite simple groups $T$ (though presumably never as large as $4/5$), and if we happen to have $c>1/2$ then $ca+b$ is maximised when $b=1$ (else it is maximised when $a=0$). 
Hence, roughly speaking, if $c$ happens to be somewhat large (for a given $T$) and we happen to be able to guarantee that $x$ can be chosen with $b$ quite small, then we might have a useful bound on the fixity of $G$ to work with. 
Bounds on $c$ can certainly be determined on a case-by-case basis from standard results about involution centralisers, but in light of the rather involved necessary and sufficient conditions for $P$ to be a maximal subgroup of $G$, it is not clear to us what can be said about the cycle structure of permutations of $\{T_1,\ldots,T_r\}$ induced by involutions in $P$. 
We therefore pose the following (somewhat vaguely worded) problem.

\begin{Prob} \label{problem:TW}
Under what conditions can a primitive permutation group of type TW and degree $d$ be guaranteed to have large fixity, where by ``large'' we mean, say, $d^{3/4}$ or more?
\end{Prob}

\appendix
\section{Additional proofs} \label{ss:app}

\noindent {\bf Proof of Proposition~\ref{Prop7/9}.}
Suppose that we are not in case~(iii). 
Then, by the final assertion of Lemma~\ref{substructureLemma}, we have $|\mathcal{P}_\theta| \le (s+1)(t+1)$, and we argue as in the proof of Theorem~\ref{subGQLemma-general}. 
We must show that we are either in case~(ii), or that $f(s,t) = ((s+1)(st+1))^{7/9} - (s+1)(t+1)$ is positive. 
We have $\tfrac{\partial f}{\partial t}(s,t) = (s+1)(7s-h(s,t))/h(s,t)$, where $h(s,t) = 9((s+1)(st+1))^{2/9}$. 
If $s \ge 13$ then (using also $2\le t\le s^2$) we have $h(s,t) \le 9 (\tfrac{14}{13} s)^{2/9} (\tfrac{27}{26} st)^{2/9} \le 9 (\tfrac{189}{169})^{2/9} s^{8/9}$, and so $7s-h(s,t) \ge s^{8/9}(7s^{1/9} - 9 (\tfrac{189}{169})^{2/9})$. 
The right-hand side is positive if and only if $s > (\tfrac{9}{7})^9(\tfrac{189}{169})^2 \approx 12.01 > 12$, and so it follows that $f(s,t) > 0$ for all $s \ge 13$, for all $s^{1/2} \le t \le s^2$. 
If $4 \le s \le 12$ then a direct calculation shows that $f(s,t) > 0$ for all $s^{1/2} \le t \le s^2$, so it remains to consider $s \in \{2,3\}$. 
If $s=2$ then, by the final paragraph of the proof of Theorem~\ref{subGQLemma-general}, either $t=2$ and every non-identity collineation of $\mathcal{Q}$ fixes at most $7 < |\mathcal{P}|^{7/9} = 15^{7/9} \approx 8.22$ points, or $t=4$ and we are in case~(ii).
Finally, if $s=3$ then $3^{1/2} \le t \le 3^2$, and a direct calculation shows that $f(3,t) > 0$ for $3^{1/2} \le t \le 7$. 
Moreover, $\mathcal{Q}$ cannot have order $(s,t)=(3,8)$ by Lemma~\ref{lemmaBasics}(iii). 
If $(s,t) = (3,9)$ then $\mathcal{Q}$ is the elliptic quadric $\mathsf{Q}^-(5,3)$ \cite[5.3.2]{FGQ}, and a {\sf FinInG}~\cite{FinInG} calculation shows that (up to conjugacy) there is a unique non-identity collineation $\theta$ fixing $40 > |\mathcal{P}|^{7/9} = 112^{7/9} \approx 39.25$ points. 
Moreover, $\mathcal{Q}_\theta$ is a generalised quadrangle of order $(3,3)$, and every other non-identity collineation of $\mathcal{Q}$ fixes at most $16 < 112^{7/9}$ points.
\hfill $\Box$ \medskip

\noindent {\bf Proof of Proposition~\ref{Prop.752}.}
Suppose that we are not in case~(ii) or (iii). 
Then $|\mathcal{P}_\theta| \le (s+1)(t+1)$ by the final assertion of Lemma~\ref{substructureLemma}. 
We show that $f(s,t) = ((s+1)(st+1))^{94/125} - (s+1)(t+1)$ is positive. 
We have $\tfrac{\partial f}{\partial t}(s,t) = (s+1)(94s-h(s,t))/h(s,t)$, where $h(s,t) = 125((s+1)(st+1))^{31/125}$. 
Let $a = 2.9701 \times 10^{15}$. 
Then $s \ge a$, so (using also $2\le t\le s^2$) we have $h(s,t) \le 125 (\tfrac{a+1}{a} s)^{31/125} (\tfrac{2a+1}{2a} st)^{31/125} \le 9 (\tfrac{(2a+1)(a+1)}{2a^2})^{31/125} s^{124/125}$, and hence $94s-h(s,t) \ge s^{124/125}(94s^{1/125} - 125 (\tfrac{(2a+1)(a+1)}{2a^2})^{31/125})$. 
The right-hand side is positive because $s \ge a > (\tfrac{125}{94})^{125}(\tfrac{(2a+1)(a+1)}{2a^2})^{31} \approx 2.97009 \times 10^{15}$, and it follows that $f(s,t) > 0$ for all $s \ge a$, for all $s^{1/2} \le t \le s^2$. 
\hfill $\Box$ \medskip

\noindent {\bf Proof of Proposition~\ref{t=s+2}.}
Since $s=t-2<t+3$, the final assertion of Lemma~\ref{substructureLemma} implies that $|\mathcal{P}_\theta| \le (s+1)(t+1) = (s+1)(s+3)$.
The result follows upon comparing this with $|\mathcal{P}| = (s+1)^3$.
\hfill $\Box$ \medskip

\noindent {\bf Proof of Corollary~\ref{diagonalPrelimk=2andr=3}.}
By Propositions~\ref{diagonalPrelim} and~\ref{diagonalPrelimk=2}, $T$ is one of the Lie type groups in the second column of Table~\ref{tab:summary4}. 
However, by Proposition~\ref{PropHC-rLeq3}, we must also have $|C_T(x)| < |T|^{5/9}$ for all $x \in T \setminus \{1\}$. 
We use this to show that $T$ cannot have type $\A_3^\varepsilon$, $\B_2=\C_2$, ${}^2\FF_4$ or $\G_2$. 

If $T \cong \G_2(q)$ then $q \ge 3$ (because $\G_2(2)$ is not simple), $|T|=q^6(q^6-1)(q^2-1)$, and we can choose $x \in T$ with $|C_T(x)| = q|\A_1(q)| = q^6(q^2-1)/\gcd(2,q-1)$ as in the proof of Lemma~\ref{centralisersExceptional}(ii). 
If $q$ is even or $q>19$ then $|C_T(x)| > |T|^{5/9}$, and if $q \in \{3,5,7,9,11,13,17,19\}$ then there is no solution of $|T|^2 = (s+1)(st+1)$ satisfying $s \ge 2$, $t \ge 2$ and properties~(ii) and~(iii) of Lemma~\ref{lemmaBasics}.
If $T \cong {}^2\FF_4(q)$ then $q=2^{2n+1}$ with $n \ge 1$ (because ${}^2\FF_4(2)$ is not simple and ${}^2\FF_4(2)'$ was treated in Proposition~\ref{diagonalPrelimk=2}), $|T| = q^{12}(q^6+1)(q^4-1)(q^3+1)(q-1)$ and, as in the proof of Lemma~\ref{centralisersExceptional}(ii), we can choose $x \in T$ with $|C_T(x)| = q^{10}|{}^2\B_2(q)| = q^{12}(q^2+1)(q-1)$. 
This yields $|C_T(x)| > |T|^{5/9}$ for all $q$. 
If $T \cong \PSp_4(q) \cong \Omega_5(q)$ then $q \ge 3$ (because $\PSp_2(2) \cong \Sym_6$ is not simple), $|T|=q^4(q^4-1)(q^2-1)/\gcd(2,q-1)$, and taking $x \in T$ with $|C_T(x)|=q^4(q^2-1)/\gcd(2,q-1)$ as in the proof of Lemma~\ref{centralisersClassical} yields $|C_T(x)| \ge |T|^{5/9}$ for all $q \ge 3$. 
Finally, if $T \cong \PSL^\varepsilon_4(q)$ (where $\text{L}^+ := \text{L}$ and $\text{L}^- := \text{U}$), then $|T| = q^6(q^4-1)(q^2-1)(q^3-\varepsilon)/\gcd(4,q-\varepsilon)$ and we can choose $x \in T$ with $|C_T(x)|=q^5|\text{GL}^\varepsilon_2(q)|/\gcd(4,q-\varepsilon) = q^6(q^2-1)(q-\varepsilon)/\gcd(4,q-\varepsilon)$ as in the proof of Lemma~\ref{centralisersClassical}. 
This yields $|C_T(x)| > |T|^{5/9}$ unless $\epsilon=+$ and $q=2$, and in this case there is no solution of $|T|^2 = (s+1)(st+1)$ satisfying $s \ge 2$, $t \ge 2$ and properties~(ii) and~(iii) of Lemma~\ref{lemmaBasics}. 
\hfill $\Box$ \medskip

\noindent {\bf Proof of Corollary~\ref{HSleftoverT} (continued).}
Now suppose that $T \cong \PSL^\varepsilon_{n+1}(q)$ with $n \in \{7,8\}$, and choose $x \in T$ as in the proof of Lemma~\ref{centralisersClassical}, so that \eqref{AnCentralisers1} holds. 
If $n=8$ then $|C_T(x)| > |T|^{7/9}$ for all $q \ge 2$, so Proposition~\ref{PropHC-rLeq3} gives a contradiction. 
If $n=7$ then $|C_T(x)| > |T|^{94/125}$ for all $q \ge 2$, so condition~(ii) must hold. 
That is, $s \le 2.9701 \times 10^{15}$, so
\begin{equation} \label{lastCorollary}
|T| = |\mathcal{P}| < (s+1)^4 \le (2.9701 \times 10^{15}+1)^4 < 7.78188 \times 10^{61}.
\end{equation} 
This implies that $q \le 9$, in which case there is no solution of $|T| = (s+1)(st+1)$ satisfying $s \ge 2$, $t \ge 2$ and properties~(ii) and~(iii) of Lemma~\ref{lemmaBasics}. 

If $T \cong \PSp_8(q)$ then, as per \eqref{CnCentralisers}, $|T| = q^{16}(q^2-1)(q^4-1)(q^6-1)(q^8-1)/\gcd(2,q-1)$ and we can choose $x \in T$ with $|C_T(x)| = q^{16}(q^2-1)(q^4-1)(q^6-1)/\gcd(2,q-1)$. 
This yields $|C_T(x)| > |T|^{94/125}$ for all $q \ge 2$, so again \eqref{lastCorollary} must hold. 
This implies that $q \le 53$, in which case there is no solution of $|T| = (s+1)(st+1)$ satisfying $s \ge 2$, $t \ge 2$ and properties~(ii) and~(iii) of Lemma~\ref{lemmaBasics}. 
Similarly, if $T \cong \Omega_9(q)$ then we can take $x \in T$ as in the proof of Lemma~\ref{centralisersClassical}, so that with $|C_T(x)| = |\SO_{2n}^\pm(q)| = q^{12}(q^4 \pm 1)(q^2-1)(q^4-1)(q^6-1)$ as in \eqref{BnCentralisers}. 
This yields $|C_T(x)| > |T|^{94/125}$ for all $q \ge 2$, so \eqref{lastCorollary} must hold, and we immediately have a contradiction because $|\Omega_9(q)| = |\PSp_8(q)|$.

Finally, suppose that $T \cong \text{P}\Omega_{2n}^\pm(q)$ with $n \in \{7,8\}$, and choose $x \in T$ as in the proof of Lemma~\ref{centralisersClassical}. 
If $n=8$ then by using \eqref{DnCentralisers1} (for $q$ odd) and \eqref{DnCentralisers3} (for $q$ even), one may check that $|C_T(x)| > |T|^{94/125}$ for all $q \ge 2$. 
Hence, \eqref{lastCorollary} must hold, and this implies that $q \in \{2,3\}$, in which case there is no solution of $|T| = (s+1)(st+1)$ satisfying $s \ge 2$, $t \ge 2$ and properties~(ii) and~(iii) of Lemma~\ref{lemmaBasics}.
Now suppose that $n=7$. 
Then \eqref{lastCorollary} holds if and only if $q \le 4$, and in this case there is no solution of $|T| = (s+1)(st+1)$ satisfying $s \ge 2$, $t \ge 2$ and properties~(ii) and~(iii) of Lemma~\ref{lemmaBasics} for $q \in \{2,4\}$. 
Therefore, we must have $q \ge 5$. 
However, in this case $|C_T(x)| > |T|^{94/125}$, so we have a contradiction. 
(To check this, note that $|C_T(x)| \ge c \cdot q^{42}(q^5-1)(q^2-1)^2(q^4-1)(q^6-1)(q^8-1)$ where $c = \tfrac{1}{2}$ or $\tfrac{1}{4}$ according as $q$ is even or odd.)
\hfill $\Box$

\end{document}